\documentclass[pdflatex,sn-mathphys-num,alpha]{sn-jnl}

\usepackage{multirow}%
\usepackage{amsmath,amssymb,amsfonts}%
\usepackage{amsthm}%
\usepackage{mathrsfs}%
\usepackage[title]{appendix}%
\usepackage{xcolor}%
\usepackage{textcomp}%
\usepackage{manyfoot}%
\usepackage{booktabs}%
\usepackage{algorithm}%
\usepackage{algorithmicx}%
\usepackage{algpseudocode}%
\usepackage{listings}%
\usepackage{cases,bbm}

\usepackage[capitalise]{cleveref}

\crefname{lemma}{Lemma}{Lemmas}
\crefname{theorem}{Theorem}{Theorems}

\crefformat{equation}{\textup{#2(#1)#3}}
\crefrangeformat{equation}{\textup{#3(#1)#4--#5(#2)#6}}
\crefmultiformat{equation}{\textup{#2(#1)#3}}{ and \textup{#2(#1)#3}}
{, \textup{#2(#1)#3}}{, and \textup{#2(#1)#3}}
\crefrangemultiformat{equation}{\textup{#3(#1)#4--#5(#2)#6}}%
{ and \textup{#3(#1)#4--#5(#2)#6}}{, \textup{#3(#1)#4--#5(#2)#6}}{, and \textup{#3(#1)#4--#5(#2)#6}}

\Crefformat{equation}{#2Equation~\textup{(#1)}#3}
\Crefrangeformat{equation}{Equations~\textup{#3(#1)#4--#5(#2)#6}}
\Crefmultiformat{equation}{Equations~\textup{#2(#1)#3}}{ and \textup{#2(#1)#3}}
{, \textup{#2(#1)#3}}{, and \textup{#2(#1)#3}}
\Crefrangemultiformat{equation}{Equations~\textup{#3(#1)#4--#5(#2)#6}}%
{ and \textup{#3(#1)#4--#5(#2)#6}}{, \textup{#3(#1)#4--#5(#2)#6}}{, and \textup{#3(#1)#4--#5(#2)#6}}

\theoremstyle{thmstyleone}%
\newtheorem{theorem}{Theorem}[section]
\newtheorem{proposition}{Proposition}[section]%

\theoremstyle{thmstyletwo}%
\newtheorem{remark}{Remark}[section]%

\theoremstyle{thmstylethree}%
\newtheorem{definition}{Definition}[section]%

\newtheorem{lemma}{Lemma}[section]

\newtheorem{hypothesis}{Hypothesis}[section]

\begin{document}

\title[Numerical Analysis of 2D Stochastic Navier--Stokes Equations with Transport Noise: Regularity and Spatial Semidiscretization
  ]{Numerical Analysis of 2D Stochastic Navier--Stokes Equations with Transport Noise: Regularity and Spatial Semidiscretization}


\author[1]{\fnm{Binjie} \sur{Li}}\email{libinjie@scu.edu.cn}
\author*[2]{\fnm{Qin} \sur{Zhou}}\email{zqmath@cwnu.edu.cn}

\affil[1]{\orgdiv{School of Mathematics}, \orgname{Sichuan University}, \city{Chengdu}, \postcode{610064}, \state{Sichuan Province}, \country{China}}
\affil[2]{\orgdiv{School of Mathematics}, \orgname{China West Normal University}, \city{Nanchong}, \postcode{637002}, \state{Sichuan Province}, \country{China}}



\abstract{\unboldmath
  This paper establishes a strong convergence rate for the spatial finite element discretization of
  a two-dimensional stochastic Navier--Stokes
  system with transport noise and no-slip boundary conditions on a convex polygonal domain.
  The main challenge arises from the lack of spatial $D(A)$-regularity of the solution
  (where $A$ is the Stokes operator),
  which prevents the direct applicability of standard error analysis techniques.
  Under a smallness assumption on the transport noise,
  for a MINI-element spatial semidiscretization, we establish a strong convergence rate
  of order $h^{\varrho} \log\big(1 + \frac{1}{h}\big)$, where $\varrho \in (0,1/2)$.
  }





\maketitle

\section{Introduction}
This paper is concerned with the numerical analysis of the two-dimensional incompressible stochastic
Navier--Stokes equations (SNSEs) driven by transport noise.
Let $0 < T < \infty$ and let $\mathcal{O} \subset \mathbb{R}^2$ be a bounded convex polygonal domain.
We consider the SNSE in the Stratonovich form:
\begin{equation}
  \label{eq:model}
  \begin{aligned}
    \mathrm{d}u(t) &= \left[ \Delta u(t) - (u(t) \cdot \nabla)u(t) - \nabla\pi(t) \right] \mathrm{d}t \\
    & \qquad {} + \sum_{n=1}^\infty \left[ (\zeta_n \cdot \nabla)u(t) - \nabla\widetilde{\pi}_n(t) \right] \circ \, \mathrm{d}W_n(t), \quad t \in [0,T],
  \end{aligned}
\end{equation}
equipped with the incompressibility constraint $\nabla \cdot u(t) = 0$, the no-slip boundary condition $u(t) = 0$ on $\partial\mathcal{O}$, and the initial condition $u(0) = u_0$.
Here, $u$ denotes the velocity field, $\pi$ the pressure, $(\zeta_n)_{n \geqslant 1}$ a sequence of divergence-free vector fields, and $(W_n)_{n \geqslant 1}$ a sequence of independent real-valued Brownian motions. The auxiliary pressure processes $(\widetilde{\pi}_n)_{n \geqslant 1}$ are introduced to ensure that the stochastic term remains divergence-free.

The systematic analysis of stochastic Navier--Stokes equations (SNSEs) was initiated in the 1970s by the seminal work of Bensoussan and Temam \cite{Bensoussan1973}.
These equations have since been extensively studied in various mathematical settings;
see, for instance, \cite{Brzezniak1991,Capinski1991,Flandoli1995martingale,Rockner2009} and the references therein.
A particularly significant class is driven by transport noise, which models the influence of unresolved,
small-scale turbulent fluctuations on the resolved, large-scale flow through a Lagrangian advection mechanism.
This formulation has a rigorous physical foundation derived from Holm's stochastic variational approach \cite{Holm2015variational},
which yields SPDEs that preserve the geometric structures of ideal fluid dynamics, such as Kelvin's circulation theorem.
The resulting Stratonovich-form SNSEs with multiplicative transport noise, for which \eqref{eq:model} serves as a canonical example,
have recently attracted significant analytical interest (see \cite{Flandoli2021high,Flandoli2023} and the references therein).

Substantial progress has also been achieved in the numerical analysis of these equations.
For problems endowed with space-periodic boundary conditions, a succession of works
\cite{Bessaih2014,Bessaih2019,Bessaih2022,Breit2021,Prohl2012,Hausenblas2019} has established strong convergence rates in probability
for various approximation schemes. The prevailing analytical strategy involves partitioning the probability space into a regular set,
where strong convergence is attained, and a singular set whose probability measure diminishes as the discretization is refined.
More recently, Feng and Vo \cite{feng2025full} derived full moment error estimates in strong norms for numerical approximations of
a two-dimensional SNSE under space-periodic boundary conditions. Their analysis relies critically on a uniform boundedness
condition for the noise coefficients, which is instrumental in proving the exponential stability of the solution.

In contrast, the analysis under no-slip boundary conditions is significantly more challenging, 
with the literature being considerably less developed.
Convergence for multiplicative noise was first shown by Brze\'zniak, Carelli, and Prohl \cite{Prohl2013} without explicit rates,
and subsequently, a first-order convergence rate was established for additive noise by Breit and Prohl \cite{Breit2023IMA}.
A more recent contribution by Breit and Prohl \cite{Breit2024} derived strong convergence rates in probability
for multiplicative noise via discrete stopping times.
This result, however, is subject to restrictive assumptions: the spatial domain must have a smooth boundary, and the noise coefficients
must be divergence-free and vanish on the boundary.
Further related results can be found in \cite{Doghman2024,Ondrejat2023}.

The challenge is further amplified when considering SNSEs driven by transport noise,
a regime that remains largely unexplored.
To our knowledge, the only convergence result is due to Breit et al.~\cite{Breit2026}, who established a mean-square rate of order $1/2$ for the temporal
discretization of the two-dimensional SNSEs under periodic boundary conditions.
However, this result relies on the restrictive assumption that the transport noise vector fields $(\zeta_n)_{n \geqslant 1}$ are constant.
This assumption is crucial for guaranteeing both the high regularity of the continuous solution and the stability of the temporal discretization,
which in turn facilitates the control of the non-Lipschitz nonlinearity.

This exposes a fundamental gap: for SNSEs with transport noise and no-slip boundaries,
solutions generally fail to possess $D(A)$-regularity,
thereby invalidating the analytical tools employed in all previous convergence analyses,
where $A$ is the Stokes operator with homogeneous Dirichlet boundary conditions.
Consequently, no convergence rates are available for spatial discretizations in this physically relevant setting.

This paper bridges this gap by introducing novel theoretical and numerical tools,
 delivering the first convergence analysis for this class of problems. Our contributions are twofold:

\begin{enumerate}
  \item First, we establish key regularity properties for the solution to the model problem \eqref{eq:model},
  assuming a smallness condition on the transport noise vector fields $(\zeta_n)_{n\geqslant 1}$,
    i.e., $\sum_{n=1}^\infty \| \zeta_n \|_{\mathbb{W}^{1,\infty}}^2$ is sufficiently small.
    This mathematical assumption has a clear physical interpretation: it characterizes a dynamical system whose energy is dominated by large-scale motions,
    which is a typical regime where the deterministic dynamics prevail and the stochastic effects can be treated as a perturbation.
    Under this condition, we demonstrate in \cref{prop:regu} that for some $\varrho \in (0,1/2)$,
    the solution $u$ satisfies
    \[
      u \in L^2(\Omega;C([0,T];X_{\varrho})) \cap L^2(\Omega \times [0,T];X_{1+\varrho}).
    \]
    This regularity result provides the essential foundation for the subsequent numerical analysis.

  \item Second, we establish a strong error estimate for a spatial semidiscretization based on
    the MINI finite element method.
    A key difficulty arises from the lack of $D(A)$-regularity of the continuous solution $u$,
    which prevents the direct application of the discrete Helmholtz projection $\mathcal{P}_h$ to $Au$.
    To overcome this obstacle, we introduce an auxiliary operator
    $\mathcal{J}_{h,\alpha} := A_h^{\alpha} \mathcal{P}_h A^{-\alpha}$ with $\alpha \in (0,1)$, where $ A_h $ is the discrete
    Stokes operator.
    This operator serves as an effective bridge between the continuous and discrete spaces, replacing the role of $\mathcal{P}_h$.
    By employing stochastic analytic tools, including It\^o's formula and the Burkholder--Davis--Gundy inequality, we derive the following strong error estimate:
    \[
      \|u - u_h\|_{L^2(\Omega; C([0,T]; L^2(\mathcal{O};\mathbb{R}^2)))}
      + \|\nabla(u - u_h)\|_{L^2(\Omega \times [0,T]; L^2(\mathcal{O};\mathbb{R}^{2\times2}))} \leqslant c h^{\varrho}\log(1+1/h),
    \]
    where $u_h$ is the strong solution of the spatial semidiscretization and $c > 0$ is a constant independent of the mesh size $h$, as stated in \cref{thm:u-uh}.
\end{enumerate}

These theoretical results provide a foundation for the further numerical analysis of fully discrete finite element approximations of the two-dimensional stochastic
Navier--Stokes equations with transport noise under no-slip boundary conditions.  
The primary novelty of this work lies in a new analytical framework that circumvents the lack of $ D(A) $-regularity—a fundamental obstacle in the numerical analysis of stochastic
Navier--Stokes equations under no-slip (i.e., homogeneous Dirichlet) boundary conditions.  
The same methodology extends naturally to the stochastic Stokes equations with Dirichlet boundary conditions,
where solution paths likewise lack $D(A)$-regularity.

The remainder of this paper is organized as follows.
Section~\ref{sec:preliminary} lays the foundation by introducing the functional--probabilistic framework,
including Sobolev and interpolation spaces, the cylindrical Brownian motion $ W $,
relevant spaces of stochastic processes, and the Stokes operator $ A $.
Section~\ref{sec:solution} establishes key regularity properties of the solution, which play an essential role in the subsequent numerical analysis.
In Section~\ref{sec:spatial}, we derive a strong convergence rate for the spatial semidiscretization based on the MINI finite element method.
The paper concludes with Section~\ref{sec:conclusions}, which summarizes the main contributions and identifies several prominent open problems for future research.

\section{Preliminaries}
\label{sec:preliminary}
\subsection{Functional framework}
Unless otherwise stated, all Banach spaces in this paper are assumed to be real.
For Banach spaces $E_1$ and $E_2$, let $\mathcal{L}(E_1,E_2)$ denote the space of
bounded linear operators from $E_1$ to $E_2$, and let $I$ denote the identity
operator on the Banach space under consideration.
Given real Banach spaces $E_1$ and $E_2$ with a continuous dense embedding $E_2 \hookrightarrow E_1$,
for each $\theta \in [0,1]$ we define
\[
  [E_1,E_2]_{\theta}
  := \bigl\{ x \in E_1 \colon x + i\,0 \in [E_{1,\mathbb{C}},E_{2,\mathbb{C}}]_{\theta} \bigr\},
\]
where $E_{1,\mathbb{C}}$ and $E_{2,\mathbb{C}}$ denote the complexifications of $E_1$
and $E_2$, respectively (see \cite[Section~B.4]{Hytonen2016}), $x + i\,0$ stands
for the canonical image of $x$ in $E_{1,\mathbb{C}}$, and
$[E_{1,\mathbb{C}},E_{2,\mathbb{C}}]_{\theta}$ denotes the corresponding complex
interpolation space (cf. \cite[Chapter~1, Section~5.1]{Yagi2010}). The space
$[E_1,E_2]_{\theta}$ is endowed with the norm inherited from
$[E_{1,\mathbb{C}}, E_{2,\mathbb{C}}]_{\theta}$, i.e.,
$ \|x\|_{[E_1,E_2]_{\theta}} := \|x + i\,0\|_{[E_{1,\mathbb{C}}, E_{2,\mathbb{C}}]_{\theta}}$,
$x \in [E_1,E_2]_{\theta}$.

Given a compact interval $D \subset \mathbb{R}$ and a Banach space $X$,
$C(D;X)$ denotes the space of continuous functions from $D$ to $X$,
endowed with the supremum norm
\[
\|f\|_{C(D;X)} = \max_{t \in D} \|f(t)\|_X, \quad f \in C(D;X).
\]

Throughout, $\mathcal{O} \subset \mathbb{R}^{2}$ is a
bounded convex polygonal domain with boundary $\partial\mathcal{O}$. For
integer $m \geqslant 1$ and $q \in [1,\infty]$, we denote by $W^{m,q}(\mathcal{O})$ and
$W^{m,q}(\mathcal{O};\mathbb{R}^{2})$ the standard scalar- and vector-valued
Sobolev spaces (cf. \cite[Chapter~III]{Boyer2013}). 
For brevity, we introduce the shorthand notation
\[
\mathbb{L}^{q} := L^{q}(\mathcal{O};\mathbb{R}^{2}),
\qquad
\mathbb{W}^{m,q} := W^{m,q}(\mathcal{O};\mathbb{R}^{2}),
\qquad m \geqslant 1, \, q \in [1,\infty].
\]

For $s \in [0,2]$, we define
\[
H^{s} := \bigl[ L^{2}(\mathcal{O}),\, W^{2,2}(\mathcal{O}) \bigr]_{s/2},
\qquad
\mathbb{H}^{s} := \bigl[ \mathbb{L}^{2},\, \mathbb{W}^{2,2} \bigr]_{s/2}.
\]
In the particular case $s = 1$, the spaces $H^{1}$ and $\mathbb{H}^{1}$ coincide,
up to equivalent norms, with $W^{1,2}(\mathcal{O})$ and
$\mathbb{W}^{1,2}$, respectively (see \cite[Theorem~1.35]{Yagi2010}).

Let $\Delta$ denote the componentwise Laplace operator on vector-valued
functions in $\mathbb{L}^2$, subject to homogeneous Dirichlet boundary
conditions, with domain
\[
  D(\Delta) := \bigl\{ v \in \mathbb{H}^2 : v|_{\partial\mathcal{O}} = 0 \bigr\}.
\]
The operator $-\Delta$ is then self-adjoint and positive definite on
$\mathbb{L}^2$. For each $\theta\in[-2,2]$, the space $\mathbb{H}_{D}^{\theta}$
is defined as the completion of $D(\Delta)$ with respect to the norm
\[
  \|v\|_{\mathbb{H}_{D}^{\theta}}
  := \bigl\| (-\Delta)^{\theta/2} v \bigr\|_{\mathbb{L}^2},
  \qquad v \in D(\Delta),
\]
where the fractional power $(-\Delta)^{\theta/2}$ is defined by the spectral theorem.
Since $\mathbb{H}_{D}^{2}$ embeds continuously into $\mathbb{H}^{2}$
(cf.\ \cite[Chapter~I, Theorem~1.8]{Girault1986}), while
$\mathbb{H}_{D}^{0} = \mathbb{L}^2 = \mathbb{H}^{0}$
and the norms coincide,
complex interpolation yields the continuous embeddings
\begin{equation}
  \label{eq:HD-to-H}
  \mathbb{H}_{D}^{\theta} \hookrightarrow \mathbb{H}^{\theta},
  \qquad \theta\in[0,2].
\end{equation}
Moreover, for $0<\theta\leqslant2$, the $\mathbb L^2$-inner product induces an isometric identification
$ \mathbb{H}_{D}^{-\theta} \cong (\mathbb{H}_{D}^{\theta})^{*} $,
where $(\mathbb{H}_{D}^{\theta})^{*}$ denotes the dual space of $\mathbb{H}_{D}^{\theta}$.

Define
\[
  \mathbb L_{\sigma}^{2}
  := \overline{\{ v \in C_{\mathrm{c}}^{\infty}(\mathcal O;\mathbb R^{2}) :
  \operatorname{div} v=0 \}}^{\,\mathbb L^{2}},
\]
i.e., the closure in $\mathbb L^{2}$ of smooth compactly supported divergence-free vector fields, and let
$\mathcal P\colon \mathbb L^{2}\to \mathbb L_{\sigma}^{2}$ denote the Helmholtz projection from $\mathbb L^{2}$
onto $\mathbb L_{\sigma}^{2}$.
The Stokes operator $A$ is defined by
\[
  A v := -\mathcal P\Delta v,
  \qquad D(A) := \mathbb L_{\sigma}^{2}\cap \mathbb{H}_{D}^{2}.
\]
The operator $A$ is self-adjoint and positive definite on $\mathbb L_{\sigma}^{2}$
(see \cite[Theorem~2.1.1]{Sohr2012navier}).

In analogy with the scale $(\mathbb{H}_{D}^{\theta})_{\theta\in[-2,2]}$,
we introduce the Hilbert scale $(X_{\theta})_{\theta\in[-2,2]}$ generated
by the Stokes operator $A$. Specifically, for each $\theta\in[-2,2]$, let
$X_{\theta}$ be the completion of $D(A)$ with respect to the norm
\[
  \|v\|_{X_{\theta}}
  := \|A^{\theta/2}v\|_{\mathbb{L}_{\sigma}^{2}},
  \qquad v \in D(A).
\]
Since $X_{2}$ embeds continuously into $\mathbb{H}_{D}^{2}$ by
\cite[Theorem~5.5]{Dauge1989}, while $X_{0} = \mathbb{L}_{\sigma}^2$ is isometrically embedded 
into $\mathbb{H}^{0} = \mathbb{L}^2$, complex interpolation yields the continuous
embeddings
\begin{equation}
  \label{eq:X-to-HD}
  X_{\theta} \hookrightarrow \mathbb{H}_{D}^{\theta},
  \qquad \theta\in[0,2].
\end{equation}
Moreover, for every $0<\theta\leqslant2$, the $\mathbb{L}_{\sigma}^{2}$-inner
product identifies $X_{-\theta}$ isometrically with the dual space of
$X_{\theta}$.

The continuous embedding $X_{1}\hookrightarrow \mathbb{H}_{D}^{1}$ implies that
$\mathcal{P}$ admits a unique extension to a bounded linear operator
$\mathcal{P}:\mathbb{H}_{D}^{-1}\to X_{-1}$, determined by the duality relation
\[
\langle \mathcal{P}v,\, w \rangle_{X_{-1},\, X_{1}}
=
\langle v,\, w \rangle_{\mathbb{H}_{D}^{-1},\, \mathbb{H}_{D}^{1}}
\qquad
\text{for all } v\in \mathbb{H}_{D}^{-1},\ w\in X_{1}.
\]
Here $\langle\cdot,\cdot\rangle_{X_{-1},X_{1}}$ denotes the duality pairing between
$X_{-1}$ and $X_{1}$, while $\langle\cdot,\cdot\rangle_{\mathbb{H}_{D}^{-1},\mathbb{H}_{D}^{1}}$
denotes that between $\mathbb{H}_{D}^{-1}$ and $\mathbb{H}_{D}^{1}$.

  Moreover, the Helmholtz projection $\mathcal{P}$ satisfies the following well-known boundedness properties.
  \begin{lemma}\label{lem:P}
    The Helmholtz projection $\mathcal{P}$ is bounded in the following senses:
    \begin{enumerate}
      \item[\textup{(i)}] $\mathcal{P} \in \mathcal{L}(\mathbb{H}^{\alpha}, \mathbb{H}^{\alpha})$ for all $\alpha \in [0,1]$.
      \item[\textup{(ii)}] $\mathcal{P} \in \mathcal{L}(\mathbb{H}^{\alpha}, X_{\beta})$ for all $\alpha \in (0,1]$ and $\beta \in (0, \alpha]$ with $\beta < 1/2$.
      \item[\textup{(iii)}] $\mathcal{P} \in \mathcal{L}(\mathbb{H}_{D}^{\alpha}, X_{\alpha})$ for all $\alpha \in [-1,0]$.
    \end{enumerate}
  \end{lemma}
  
  \begin{proof}
    \textbf{Proof of (i).}
    Since $\mathcal{P}$ is the orthogonal projection onto $\mathbb{L}_{\sigma}^{2}$,
    it is trivially bounded on $\mathbb{H}^{0} = \mathbb{L}^2$.
    Moreover, on convex polygonal domains, the elliptic regularity theory for the
    Poisson equation with Neumann boundary conditions
    (cf.~\cite[Chapter~I, Theorem~1.10]{Girault1986}) yields the
    $\mathbb{H}^{1}$-boundedness of $\mathcal{P}$ via the argument in
    \cite[Proposition~IV.3.7]{Boyer2013}.
    Property~(i) follows by complex interpolation between these two endpoint
    estimates.
  
    \textbf{Proof of (ii).}
    This is a direct consequence of \cite[Proposition~2.16 and Theorem~5.1]{Mitrea2008}.
  
    \textbf{Proof of (iii).}
    The continuous embedding
    $X_{1} \hookrightarrow \mathbb{H}_{D}^{1}$
    implies that the extension of $\mathcal{P}$ defined above is bounded from
    $\mathbb{H}_{D}^{-1}$ to $X_{-1}$, i.e.,
    $\mathcal{P} \in \mathcal{L}(\mathbb{H}_{D}^{-1}, X_{-1})$.
    Moreover, since $\mathcal{P}$ is the $\mathbb{L}^{2}$-orthogonal projection
    onto $\mathbb{L}_{\sigma}^{2}$ and
    $\mathbb{H}_{D}^{0} = \mathbb{L}^{2}$ and
    $X_{0} = \mathbb{L}_{\sigma}^{2}$ hold with identical norms,
    we have $\mathcal{P} \in \mathcal{L}(\mathbb{H}_{D}^{0}, X_{0})$
    with operator norm one.
    For $\alpha \in [-1,0]$, complex interpolation with parameter $1+\alpha$
    between these two endpoint estimates, using the isometric identifications
    $[\mathbb{H}_{D}^{-1}, \mathbb{H}_{D}^{0}]_{1+\alpha} = \mathbb{H}_{D}^{\alpha}$
    and
    $[X_{-1}, X_{0}]_{1+\alpha} = X_{\alpha}$, yields property~(iii).
  \end{proof}

\subsection{Probabilistic framework}
Fix $T > 0$. Let $(\Omega, \mathcal{F}, \mathbb{P})$ be a probability space,
endowed with a filtration $\mathbb{F} = (\mathcal{F}_t)_{0 \leqslant t \leqslant T}$
satisfying the usual conditions (i.e., $\mathbb{F}$ is right-continuous and
$\mathcal{F}_0$ contains all $\mathbb{P}$-null sets of $\mathcal{F}$),
and supporting a sequence $(W_n)_{n \geqslant 1}$ of independent real-valued
standard $\mathbb{F}$-Brownian motions. Let $\ell^2$ denote the real Hilbert space
of square-summable sequences. We define the $\ell^2$-cylindrical Brownian motion
$W = (W(t))_{0 \leqslant t \leqslant T}$ by
\[
W(t)l := \sum_{n = 1}^{\infty} l_n\, W_n(t), \qquad t \in [0,T],
\; l = (l_n)_{n \geqslant 1} \in \ell^2.
\]
For each $t \in [0,T]$ and every $l \in \ell^2$,
the series converges in $L^2(\Omega, \mathcal{F}_t, \mathbb{P})$,
so that $W(t)$ is a well-defined bounded linear operator from $\ell^2$ into
$L^2(\Omega, \mathcal{F}_t, \mathbb{P})$ with operator norm $\sqrt{t}$;
in particular, $W(t)l$ is $\mathcal{F}_t$-measurable and centered Gaussian with
variance $t\|l\|_{\ell^2}^2$.

For a real separable Hilbert space $U$, we denote by $\mathcal{L}_2(\ell^2, U)$ the space of
Hilbert--Schmidt operators from $\ell^2$ to $U$ (cf.~\cite[Appendix B]{Liu2015}).

Let $E$ be a separable Banach space. For $p \in [1, \infty]$, we denote by $L^p_{\mathcal{F}_0}(\Omega; E)$ and $L^p(\Omega; E)$ the Bochner spaces $L^p(\Omega, \mathcal{F}_0, \mathbb{P}; E)$ and $L^p(\Omega, \mathcal{F}, \mathbb{P}; E)$, respectively.
For $p \in [1, \infty)$, we write $L^p(\Omega \times [0,T]; E)$ for the Bochner space
\[
L^p\bigl(\Omega \times [0,T],\ \mathcal{F}_T \otimes \mathcal{B}([0,T]),\ \mathbb{P} \otimes \lambda;\ E\bigr),
\]
where $\lambda$ denotes Lebesgue measure on $[0,T]$.
For $p \in [1,\infty)$, we define $L^p_{\mathbb F}(\Omega\times[0,T];E)$ as the set of
equivalence classes in $L^p(\Omega\times[0,T];E)$ that admit an $\mathbb F$-progressively
measurable representative. This set is a closed subspace of $L^p(\Omega\times[0,T];E)$.

\section{Well-posedness and regularity}
\label{sec:solution}
The stochastic model \eqref{eq:model} is originally formulated in the Stratonovich form,
which arises naturally from its physical derivation involving
transport noise. Applying the Helmholtz projection $\mathcal{P}$ yields the following
abstract stochastic evolution equation in the solenoidal space $\mathbb{L}_{\sigma}^2$:
\[
  \mathrm{d}u(t) = - \bigl( A u(t) + \mathcal{P}((u(t) \cdot \nabla)u(t)) \bigr) \,\mathrm{d}t
  + \sum_{n=1}^\infty L_{\zeta_n}u(t) \circ \, \mathrm{d}W_n(t), \quad t \in [0,T],
\]
where the transport operator $ L_{\zeta_n} $ is formally defined by
\begin{equation}
  \label{eq:L-zeta-n}
  L_{\zeta_n} v := \mathcal{P}\bigl((\zeta_n \cdot \nabla) v \bigr).
\end{equation}
To enable a rigorous analytical treatment, we convert this equation to its It\^o form.
This transformation introduces an It\^o correction term in the drift:
\begin{align*}
  \mathrm{d}u(t) &= \Bigl(- A u(t) - \mathcal{P}((u(t)\cdot\nabla)u(t))
    + \frac12\sum_{n=1}^\infty  L_{\zeta_n}^2 u(t)
  \Bigr) \,\mathrm{d}t + \sum_{n=1}^{\infty} L_{\zeta_n}u(t) \, \mathrm{d}W_n(t).
\end{align*}

This additional term motivates the introduction of the modified Stokes operator
\begin{equation}
  \label{eq:calA-def}
  \mathcal{A}v := A v - \frac{1}{2} \sum_{n=1}^\infty L_{\zeta_n}^2v,
  \quad \forall v \in X_{2}.
\end{equation}
With this definition, the It\^o formulation simplifies to the compact form:
\begin{equation}
    \label{eq:Ito-form}
  \mathrm{d}u(t) = \bigl(- \mathcal{A} u(t) - \mathcal{P}((u(t)\cdot\nabla)u(t)) \bigr) \,\mathrm{d}t
  + \sum_{n=1}^{\infty} L_{\zeta_n} u(t) \,\mathrm{d}W_n(t), \quad t \in [0,T].
\end{equation}

Throughout this paper, we impose the following structural and smallness assumptions on the transport noise vector
fields $(\zeta_n)_{n \geqslant 1}$. 
\begin{hypothesis}
  \label{hyp:zeta}
  Let $(\zeta_n)_{n \geqslant 1} \subset \mathbb{W}^{1,\infty}$ be divergence-free vector fields satisfying
  the summability condition
  \begin{equation}
    \label{eq:zeta-sum-cond} 
    \sum_{n=1}^\infty \|\zeta_n\|_{\mathbb{W}^{1,\infty}}^2 < \infty. 
  \end{equation}
  Moreover, we assume that the sum in \eqref{eq:zeta-sum-cond}
  is sufficiently small that
  there exists a constant $\kappa_{\zeta} \in (0,1)$ for which
  \begin{equation}
    \label{eq:sigma-key}
    \Big\| \frac{1}{2} \sum_{n=1}^\infty L_{\zeta_n}^2 v \Big\|_{\mathbb{L}_{\sigma}^2}
    \leqslant \kappa_{\zeta} \|v\|_{X_{2}}, \quad \forall v \in X_{2}.
  \end{equation}
\end{hypothesis}

\begin{remark}
  For any $\zeta \in \mathbb{W}^{1,\infty}$ and $v \in X_2$,
  the continuous embedding $X_2 \hookrightarrow \mathbb{H}^2$
  (see \cref{eq:X-to-HD,eq:HD-to-H} with $\theta=2$), together with
  \cref{lem:P}(i) applied with $\alpha=1$, yields
  \[
    \|L_{\zeta}^2 v\|_{\mathbb{L}_{\sigma}^2}
    = \big\|\mathcal{P}\big[(\zeta\cdot\nabla)
      \mathcal{P}((\zeta\cdot\nabla)v)\big]\big\|_{\mathbb{L}_{\sigma}^2}
    \leqslant c\,\|\zeta\|_{\mathbb{W}^{1,\infty}}^2 \|v\|_{X_2},
  \]
  where the constant $c>0$ is independent of $\zeta$ and $v$.
  Consequently, under condition \eqref{eq:zeta-sum-cond}, the series
  $\sum_{n=1}^\infty L_{\zeta_n}^2v$ converges absolutely in
  $\mathbb{L}_{\sigma}^2$ for every $v \in X_2$. Furthermore,
  \[
    \Big\|\frac{1}{2}\sum_{n=1}^\infty L_{\zeta_n}^2v\Big\|_{\mathbb{L}_{\sigma}^2}
    \leqslant \frac{1}{2}\sum_{n=1}^\infty
      \|L_{\zeta_n}^2v\|_{\mathbb{L}_{\sigma}^2} 
    \leqslant \frac{c}{2}
      \Big(\sum_{n=1}^\infty\|\zeta_n\|_{\mathbb{W}^{1,\infty}}^2\Big)
      \|v\|_{X_2}.
  \]
  Hence, provided the sum in \eqref{eq:zeta-sum-cond} is sufficiently small,
  inequality \eqref{eq:sigma-key} holds for some constant
  $\kappa_{\zeta}\in(0,1)$.
\end{remark}

\begin{lemma}
  Suppose that \cref{hyp:zeta} holds. Then, for every
  $\theta \in [0,1]$, there exists a constant $c>0$, depending only on
  $\theta$, $\mathcal{O}$, and $\sum_{n=1}^\infty \|\zeta_n\|_{\mathbb{W}^{1,\infty}}^2$, such that
  \begin{equation}
    \label{eq:L-zeta-n-theta-1-sum}
    \Bigl\|\sum_{n=1}^{\infty}L_{\zeta_n}^2v\Bigr\|_{X_{\theta-1}}
    \leqslant c\|v\|_{X_{\theta+1}},
    \qquad \forall v \in X_{\theta+1}.
  \end{equation}
\end{lemma}

\begin{proof}
  Throughout the proof, $c$ denotes a positive constant depending only on
  $\theta$, $\mathcal{O}$, and $\sum_{n=1}^\infty \|\zeta_n\|_{\mathbb{W}^{1,\infty}}^2$,
  whose value may change from line to line.
  
  For each $n\geqslant 1$, the following standard estimates hold:
  \begin{align*}
    \|(\zeta_n\cdot\nabla)v\|_{\mathbb{L}^2}
    &\leqslant c\|\zeta_n\|_{\mathbb{L}^\infty}\|v\|_{\mathbb{H}^1},
    \qquad v\in \mathbb{H}^1, \\
    \|(\zeta_n\cdot\nabla)v\|_{\mathbb{H}_D^{-1}}
    &\leqslant c\|\zeta_n\|_{\mathbb{L}^\infty}\|v\|_{\mathbb{L}^2},
    \qquad v\in \mathbb{L}^2.
  \end{align*}
  Interpolating between these two bounds for $\theta\in(0,1)$,
  and using the first and second estimates directly to treat the endpoint cases $\theta=1$ and $\theta=0$,
  respectively, we obtain
  \begin{equation}
    \label{eq:zetan-cdot-nabla-v}
    \|(\zeta_n\cdot\nabla)v\|_{\mathbb{H}_D^{\theta-1}}
    \leqslant c \|\zeta_n\|_{\mathbb{L}^\infty}\|v\|_{\mathbb{H}^\theta},
    \qquad v\in \mathbb{H}^\theta.
  \end{equation}
  Combining \cref{eq:zetan-cdot-nabla-v} with \cref{lem:P}(iii), applied
  with $\alpha=\theta-1$, yields
  \begin{equation}
    \label{eq:L-zeta-n-theta-1}
    \|L_{\zeta_n}v\|_{X_{\theta-1}}
    \leqslant
    c\|\zeta_n\|_{\mathbb{L}^\infty}\|v\|_{\mathbb{H}^\theta},
    \qquad v\in \mathbb{H}^\theta.
  \end{equation}
  Similarly, combining the product estimate (see \cite[Chapter~4]{Runst1996sobolev})
  \begin{equation}
    \label{eq:zetan-nabla-v-H}
    \|(\zeta_n\cdot\nabla)v\|_{\mathbb{H}^{\theta}}
    \leqslant
    c\|\zeta_n\|_{\mathbb{W}^{1,\infty}}\|v\|_{\mathbb{H}^{\theta+1}},
    \qquad v\in X_{\theta+1},
  \end{equation}
  with the continuous embedding $X_{\theta+1} \hookrightarrow \mathbb{H}^{\theta+1}$
  (see \cref{eq:X-to-HD,eq:HD-to-H} with $\theta$ replaced by $\theta+1$)
  and \cref{lem:P}(i) (with $\alpha=\theta$) gives
  \begin{equation}
    \label{eq:L-zeta-n-theta}
    \|L_{\zeta_n}v\|_{\mathbb{H}^\theta}
    \leqslant
    c\|\zeta_n\|_{\mathbb{W}^{1,\infty}}\|v\|_{X_{\theta+1}},
    \qquad v\in X_{\theta+1}.
  \end{equation}

  Now let $v\in X_{\theta+1}$. Applying
  \cref{eq:L-zeta-n-theta-1,eq:L-zeta-n-theta} successively and using
  $\|\zeta_n\|_{\mathbb{L}^\infty}
  \leqslant\|\zeta_n\|_{\mathbb{W}^{1,\infty}}$, we obtain
  \begin{align*}
    \|L_{\zeta_n}^2v\|_{X_{\theta-1}}
    &\leqslant
    c\|\zeta_n\|_{\mathbb{L}^\infty}
    \|L_{\zeta_n}v\|_{\mathbb{H}^\theta} \\
    &\leqslant
    c\|\zeta_n\|_{\mathbb{L}^\infty}
    \|\zeta_n\|_{\mathbb{W}^{1,\infty}}
    \|v\|_{X_{\theta+1}} \\
    &\leqslant
    c\|\zeta_n\|_{\mathbb{W}^{1,\infty}}^2
    \|v\|_{X_{\theta+1}}.
  \end{align*}
  Summing over $n$ and invoking \cref{eq:zeta-sum-cond}, we deduce that the
  series $\sum_{n=1}^{\infty}L_{\zeta_n}^2v$ converges absolutely in
  $X_{\theta-1}$ for every $v\in X_{\theta+1}$ and satisfies
  \cref{eq:L-zeta-n-theta-1-sum}. This completes the proof.
\end{proof}

Under \cref{hyp:zeta}, the operator $\mathcal{A}$ is self-adjoint and positive definite on
$\mathbb{L}_{\sigma}^{2}$, with domain $D(\mathcal{A})=X_{2}$.
Moreover, it is an isomorphism from $X_{2}$ onto $\mathbb{L}_{\sigma}^{2}$
and satisfies the following two-sided estimate:
\begin{equation}
  \label{eq:calH2-dotH2-equiv}
  (1-\kappa_{\zeta})\|v\|_{X_{2}}
  \leqslant \|\mathcal{A}v\|_{\mathbb{L}_{\sigma}^{2}}
  \leqslant (1+\kappa_{\zeta})\|v\|_{X_{2}},
  \qquad v\in X_{2}.
\end{equation}
Analogously to the scale $(X_{\theta})_{\theta\in[-2,2]}$, let
$(\mathcal{X}_{\theta})_{\theta\in[-2,2]}$ denote the scale generated by
$\mathcal{A}$. By \eqref{eq:calH2-dotH2-equiv} and complex interpolation
(see, e.g., \cite[Theorems~1.15 and 16.1]{Yagi2010}),
$X_\theta=\mathcal{X}_\theta$ for every $\theta\in[0,2]$, and the norms
$\|\cdot\|_{\mathcal{X}_{\theta}}$ and $\|\cdot\|_{X_{\theta}}$ are
equivalent; more precisely,
\begin{equation}
  \label{eq:X-calX-equiv}
  (1-\kappa_{\zeta})^{\theta/2}\|\cdot\|_{X_{\theta}}
  \leqslant \|\cdot\|_{\mathcal{X}_{\theta}}
  \leqslant (1+\kappa_{\zeta})^{\theta/2}\|\cdot\|_{X_{\theta}},
  \qquad \theta \in [0,2].
\end{equation}
By duality, for each $\theta \in [-2,0)$, $X_{\theta}$ coincides with
$\mathcal{X}_{\theta}$, and
\begin{equation}
  \label{eq:X-calX-equiv-negative}
  (1+\kappa_{\zeta})^{\theta/2}\| \cdot \|_{X_{\theta}}
  \leqslant \| \cdot \|_{\mathcal{X}_{\theta}}
  \leqslant (1-\kappa_{\zeta})^{\theta/2}\|\cdot\|_{X_{\theta}}.
\end{equation}
Moreover, \cref{eq:X-calX-equiv,eq:X-to-HD,eq:HD-to-H} yield the
continuous embeddings
\begin{equation} \label{eq:calX-to-H}
  \mathcal{X}_{\theta} \hookrightarrow \mathbb{H}^{\theta},
  \quad \theta \in [0,2].
\end{equation}

For $v\in\mathcal{X}_{2}=X_{2}$, integration by parts—together with the divergence-free condition
on $\zeta_{n}$, the boundary condition $v|_{\partial\mathcal{O}}=0$, and the fact that $\mathcal{P}$
denotes the $\mathbb{L}^{2}$-orthogonal projection onto $\mathbb{L}_{\sigma}^{2}$—yields
$\langle L_{\zeta_n}^{2}v,v\rangle =-\|L_{\zeta_n}v\|_{\mathbb{L}_{\sigma}^{2}}^{2}$.
Consequently, since $\langle Av,v\rangle=\|v\|_{X_{1}}^{2}$,
\[
  \|v\|_{\mathcal{X}_{1}}^{2}
  =\langle\mathcal{A}v,v\rangle
  =\|v\|_{X_{1}}^{2}
  +\frac{1}{2}\sum_{n=1}^{\infty}
  \|L_{\zeta_n}v\|_{\mathbb{L}_{\sigma}^{2}}^{2},
\]
where $\langle\cdot,\cdot\rangle$ denotes the inner product of $\mathbb{L}_{\sigma}^{2}$.
By \cref{eq:X-calX-equiv} with $\theta=1$, the first term on the right-hand side defines
a continuous functional on $\mathcal{X}_{1}$. Likewise, the second term defines a continuous
functional on $\mathcal{X}_{1}$, as follows from \cref{eq:L-zeta-n-theta-1} with $\theta=1$,
the summability condition \cref{eq:zeta-sum-cond}, and the continuous embedding
\cref{eq:calX-to-H} with $\theta=1$. Since $\mathcal{X}_{2}$ is densely and continuously
embedded in $\mathcal{X}_{1}$, the identity therefore extends by continuity to all $v\in\mathcal{X}_{1}$, giving
\begin{equation}
  \label{eq:calH12}
  \|v\|_{\mathcal{X}_{1}}^{2}
  =\|v\|_{X_{1}}^{2}
  +\frac{1}{2}\sum_{n=1}^{\infty}
    \|L_{\zeta_n}v\|_{\mathbb{L}_{\sigma}^{2}}^{2},
  \qquad v\in\mathcal{X}_{1}.
\end{equation}

For $v\in\mathcal{X}_{1}$, define
$F(v):\ell^{2}\to\mathbb{L}_{\sigma}^{2}$ by
\begin{equation}
  \label{eq:F-def}
  F(v)l=\sum_{n=1}^{\infty}l_nL_{\zeta_n}v,
  \qquad
  l=(l_n)_{n=1}^{\infty}\in\ell^{2},
\end{equation}
where the series converges in $\mathbb{L}_{\sigma}^{2}$
(a direct consequence of \cref{eq:zeta-sum-cond},
\cref{eq:L-zeta-n-theta-1} and \cref{eq:calX-to-H} with $\theta=1$).
It then follows from \cref{eq:X-calX-equiv} with $\theta=1$ and
\cref{eq:calH12} that $F(v)$ is Hilbert--Schmidt for every
$v\in\mathcal{X}_{1}$, with
\begin{align*}
  \|F(v)\|_{\mathcal{L}_{2}(\ell^{2},\mathbb{L}_{\sigma}^{2})}^{2}
  = \sum_{n=1}^{\infty}\|L_{\zeta_n}v\|_{\mathbb{L}_{\sigma}^{2}}^{2}
  =2\|v\|_{\mathcal{X}_{1}}^{2}-2\|v\|_{X_{1}}^{2}
  \leqslant\frac{2\kappa_{\zeta}}{1+\kappa_{\zeta}}
  \|v\|_{\mathcal{X}_{1}}^{2}.
\end{align*}
Hence, $F\in
\mathcal{L}\bigl(\mathcal{X}_{1},
\mathcal{L}_{2}(\ell^{2},\mathbb{L}_{\sigma}^{2})\bigr)$, and
\begin{equation}
  \label{eq:tmp0}
  \|F\|_{\mathcal{L}(\mathcal{X}_{1},
  \mathcal{L}_{2}(\ell^{2},\mathbb{L}_{\sigma}^{2}))}
  \leqslant
  \sqrt{\frac{2\kappa_{\zeta}}{1+\kappa_{\zeta}}}.
\end{equation}

We are now in a position to define the notion of a solution.

\begin{definition}\label{def:sol-def}
  Let $u_0 \in L_{\mathcal{F}_0}^4(\Omega; \mathbb{L}_{\sigma}^2)$.
  An $\mathbb{F}$-adapted process
  $u \colon \Omega \times [0,T] \to \mathbb{L}_{\sigma}^2$
  is called a \emph{solution} to the model problem \eqref{eq:model}
  if the following conditions are satisfied:
  \begin{enumerate}[label=(\roman*)]
    \item[(i)] \label{item:paths}
      $u$ has $\mathbb{P}$-a.s.\ continuous paths in $\mathbb{L}_{\sigma}^2$;
    \item[(ii)] 
      the $\mathbb{P}\otimes\mathrm{d}t$-equivalence class of $u$ belongs to
      $L_{\mathbb{F}}^2(\Omega\times[0,T]; \mathcal{X}_1)$;
    \item[(iii)] 
      $\mathbb{P}$-a.s., the identity
      \begin{equation}\label{eq:u-def}
        u(t) = u_0
        - \int_0^t \Bigl( \mathcal{A}u(s)
          + \mathcal{P}\bigl((u(s)\cdot\nabla)u(s)\bigr) \Bigr) \,\mathrm{d}s
        + \int_0^t F(u(s)) \,\mathrm{d}W(s)
      \end{equation}
      holds in $\mathcal{X}_{-1}$ for all $t \in [0,T]$, where the noise
      coefficient $F$ is given by \eqref{eq:F-def}.
  \end{enumerate}
  Two solutions are identified whenever they are indistinguishable.
\end{definition}

The main result of this section is the following proposition.

\begin{proposition}\label{prop:regu}
  Suppose that Hypothesis \ref{hyp:zeta} holds and that
  \[
    u_0 \in L_{\mathcal{F}_0}^{\infty}(\Omega;\mathbb{L}^{2}_{\sigma}).
  \]
  Then the model problem \eqref{eq:model} admits a unique solution $u$, and
  \begin{align}
    u &\in L^{\infty}(\Omega;L^{2}(0,T;X_{1}))
    \cap L^{\infty}(\Omega;C([0,T];\mathbb{L}^{2}_{\sigma})).
    \label{eq:u-basic}
  \end{align}
  Moreover, if, in addition, $u_0 \in L_{\mathcal{F}_0}^{2}\bigl(\Omega;X_{1/2}\bigr)$,
  then there exists a constant $\varrho\in(0,\tfrac{1}{2})$ such that
  the solution $u$ enjoys the additional regularity
  \begin{align}
    u &\in L^{2}\bigl(\Omega;C([0,T];X_{\varrho})\bigr),
    \label{eq:u-C}\\
    u &\in L^{2}\bigl(\Omega\times[0,T];X_{\varrho+1}\bigr).
    \label{eq:u-L2}
  \end{align}
\end{proposition}

\begin{remark}
  We note that if the transport noise vector fields $(\zeta_n)_{n\geqslant 1} \subset \mathbb{W}^{1,\infty}$
  are divergence-free and satisfy the summability condition
  \[
    \sum_{n=1}^\infty \|\zeta_n\|_{\mathbb{L}^\infty}^2 < \infty,
  \]
  then the existence and uniqueness of a solution can be established under the significantly
  weaker assumption $u_0 \in L_{\mathcal{F}_0}^4(\Omega;\mathbb{L}_{\sigma}^2)$ (cf. \cite[Theorem~4.21]{Flandoli2023}).
  In the present work, the smallness condition on the transport noise vector fields imposed in
  \cref{hyp:zeta} (i.e., $\sum_{n=1}^\infty \|\zeta_n\|_{\mathbb{W}^{1,\infty}}^2$ being sufficiently small)
  is not required for mere existence and uniqueness.
  Instead, the smallness condition plays a crucial role in our proof
  of the higher spatial regularity $u \in L^2_{\mathbb{F}}(\Omega \times [0,T]; X_{\varrho+1})$
  as stated in \eqref{eq:u-L2}.
\end{remark}

\begin{remark}
  \label{rem:regularity}
  A fundamental obstruction to higher regularity stems from the mapping
  properties of the Helmholtz projection $\mathcal{P}$.
  Indeed, by \cref{lem:P}(ii), $\mathcal{P}$ fails to be bounded from
  $\mathbb{H}^{\alpha}$ to $X_{\beta}$ for $\beta \geqslant 1/2$;
  consequently, $F(v)$ generally cannot be expected to take values in
  $\mathcal{L}_2(\ell^2, X_{\theta})$ for any $\theta \geqslant 1/2$,
  no matter how regular $v$ may be.
  This constitutes a fundamental obstruction to attaining higher
  regularity for the solution $u$ of the model problem \cref{eq:model}.
\end{remark}

\begin{remark}
The well-posedness theory for the two-dimensional stochastic Navier--Stokes equations with transport noise is well established.  
For the Cauchy problem, Mikulevicius and Rozovskii~\cite{Mikulevicius2005} proved the existence and pathwise uniqueness of a global strong solution.  
In the periodic setting, Flandoli and Luongo~\cite[Theorem 4.21]{Flandoli2023} established the existence and uniqueness of weak
solutions, while Breit et al.~\cite{Breit2026} derived certain regularity estimates for weak pathwise solutions under the
assumption that the transport noise vector fields $(\zeta_n)_{n \geqslant 1}$ are constant.  
The analysis in~\cite[Theorem 4.21]{Flandoli2023} also extends to no-slip boundary conditions with only minor modifications.  
However, to the best of our knowledge, the regularity results \cref{eq:u-C,eq:u-L2} in \cref{prop:regu} are new,
which are crucial for our numerical analysis.
We also refer the reader to the recent work of Agresti and Veraar \cite{Agresti2024} for the theoretical
analysis of a class of stochastic Navier--Stokes equations that are closely related to the model problem \cref{eq:model}.
\end{remark}

\begin{proof}[\textbf{Proof of \cref{prop:regu}}]
The existence and uniqueness of a solution, along with the basic regularity \eqref{eq:u-basic}, are standard.
These results follow directly from a variational approach, mirroring the proof of \cite[Theorem~4.21]{Flandoli2023}
with minor adjustments to accommodate our specific setting.
Our contribution is the derivation of the higher-order regularity estimates \eqref{eq:u-C} and \eqref{eq:u-L2}.
Throughout the remainder of the proof, $c > 0$ denotes a generic constant,
independent of the indices $k,n \geqslant 1$ and of the parameter $\epsilon > 0$,
whose value may change from one occurrence to another.

\medskip\noindent
\textbf{Step 1.} We prove the existence of $ 0 < \varrho^{*} < \tfrac{1}{2} $ such that
\begin{equation}
  \label{eq:PF}
  \sup_{1 \leqslant \alpha \leqslant \varrho^{*} + 1}
  \|F\|_{\mathcal{L}(\mathcal{X}_{\alpha},\, \mathcal{L}_{2}(\ell^2,\mathcal{X}_{\alpha-1}))} < \sqrt{2}.
\end{equation}
The argument is divided into two parts.

\textbf{(a) Bound for $\alpha = \alpha_0 > 1$.}  
Fix any $ \alpha_0 \in (1, \tfrac{3}{2}) $. For any $ n \geqslant 1 $ and $ v \in X_{\alpha_0} $,  
by Lemma~\ref{lem:P}(ii) (with $\alpha=\beta = \alpha_0-1$),
the continuous embedding $X_{\alpha_0} \hookrightarrow \mathbb{H}^{\alpha_0}$
(see \cref{eq:HD-to-H,eq:X-to-HD} with $\theta=\alpha_0$),
and \cref{eq:zetan-nabla-v-H} with $\theta=\alpha_0-1$
we have
\[
\|\mathcal{P}((\zeta_n \cdot \nabla)v)\|_{X_{\alpha_0-1}} \leqslant c \|\zeta_n\|_{\mathbb{W}^{1,\infty}} \|v\|_{X_{\alpha_0}}.
\]
Therefore,
\[
\begin{aligned}
  \|F(v)\|_{\mathcal{L}_{2}(\ell^2, X_{\alpha_0-1})}
  = \left( \sum_{n=1}^\infty \|\mathcal{P}((\zeta_n \cdot \nabla)v)\|_{X_{\alpha_0-1}}^2 \right)^{1/2} 
  \leqslant \left( c\sum_{n=1}^\infty \|\zeta_n\|_{\mathbb{W}^{1,\infty}}^2 \right)^{1/2} \|v\|_{X_{\alpha_0}}.
\end{aligned}
\]
Using the norm equivalence \eqref{eq:X-calX-equiv} with $\theta = \alpha_0-1$
and $\theta = \alpha_0$, we conclude that
\begin{equation}
  \label{eq:tmp2}
  \|F\|_{\mathcal{L}(\mathcal{X}_{\alpha_0}, \, \mathcal{L}_{2}(\ell^2, \mathcal{X}_{\alpha_0-1}))}
  \leqslant \left( c\sum_{n=1}^\infty \|\zeta_n\|_{\mathbb{W}^{1,\infty}}^2 \right)^{1/2}.
\end{equation}

\textbf{(b) Interpolation.}
We now interpolate between the bounds \eqref{eq:tmp0} and \eqref{eq:tmp2}
to derive estimates at intermediate regularity levels.
For $\theta\in(0,1)$, set $\alpha=1+(\alpha_0-1)\theta$.
Standard interpolation theory
(see, e.g., \cite[Theorem~16.1]{Yagi2010}) yields
\[
  [\mathcal{X}_{1},\mathcal{X}_{\alpha_0}]_{\theta}
  = \mathcal{X}_{\alpha},
  \qquad
  [\mathcal{L}_{2}(\ell^2,\mathcal{X}_{0}),
   \mathcal{L}_{2}(\ell^2,\mathcal{X}_{\alpha_0-1})]_{\theta}
  = \mathcal{L}_{2}(\ell^2,\mathcal{X}_{\alpha-1}),
\]
with equality of norms.
Applying the interpolation theorem for bounded linear operators
(see \cite[Theorem~1.15]{Yagi2010})
to \eqref{eq:tmp0} and \eqref{eq:tmp2}, we obtain
\[
  \|F\|_{\mathcal{L}(\mathcal{X}_{\alpha},
    \mathcal{L}_{2}(\ell^2,\mathcal{X}_{\alpha-1}))}
  \leqslant
  \|F\|_{\mathcal{L}(\mathcal{X}_{1},
    \mathcal{L}_{2}(\ell^2,\mathcal{X}_{0}))}^{1-\theta}
  \|F\|_{\mathcal{L}(\mathcal{X}_{\alpha_0},
    \mathcal{L}_{2}(\ell^2,\mathcal{X}_{\alpha_0-1}))}^{\theta}.
\]
Inserting \eqref{eq:tmp0} and \eqref{eq:tmp2} into this estimate yields
\[
  \|F\|_{\mathcal{L}(\mathcal{X}_{\alpha},
    \mathcal{L}_{2}(\ell^2,\mathcal{X}_{\alpha-1}))}
  \leqslant
  \left(\frac{2\kappa_\zeta}{1+\kappa_\zeta}\right)^{(1-\theta)/2}
  \left(
    c\sum_{n=1}^{\infty}
    \|\zeta_n\|_{\mathbb{W}^{1,\infty}}^{2}
  \right)^{\theta/2},
\]
where the constant $c>0$ depends only on $\alpha_0$, $\kappa_\zeta$, and $\mathcal{O}$;
in particular, it is independent of $\theta$.
The right-hand side extends continuously to $\theta=0$, with value
$\sqrt{2\kappa_\zeta/(1+\kappa_\zeta)}<\sqrt{2}$.
Hence, by continuity, there exists $\theta^{*}\in(0,1)$ such that
the right-hand side is strictly less than $\sqrt{2}$
for every $\theta\in[0,\theta^{*}]$.
Setting
$\varrho^{*}:=(\alpha_0-1)\theta^{*}\in(0,\tfrac{1}{2})$
then completes the proof of \eqref{eq:PF}.

\medskip\noindent
\textbf{Step 2.}  
It is standard that the space $\mathbb{L}_\sigma^2$ admits an orthonormal basis $(e_n)_{n\geqslant 1}$ of
eigenvectors of the operator $\mathcal A$. For each $k\geqslant 1$, let $\mathcal P_k$ denote the
orthogonal projection in $\mathbb{L}^2$ onto the finite-dimensional subspace $\operatorname{span}\{e_1,\dots,e_k\}$.
We consider the Galerkin approximation
\[
\begin{cases}
\mathrm{d}u^{(k)}(t) =
-\mathcal P_k\big( \mathcal A u^{(k)}(t) + (u^{(k)}(t)\cdot\nabla)u^{(k)}(t) \big) \, \mathrm{d}t
+ \displaystyle\sum_{n=1}^k \mathcal P_k L_{\zeta_n}u^{(k)}(t) \, \mathrm{d}W_n(t), \quad t \in [0,T], \\
u^{(k)}(0) = \mathcal P_k u_0.
\end{cases}
\]
Since this is a finite-dimensional stochastic differential equation with a locally monotone drift and linear
diffusion coefficients, the classical theory (see, e.g., \cite[Theorem~3.27]{Pardoux2014})
guarantees the existence of a unique process
$u^{(k)}\colon\Omega\times[0,T]\to \operatorname{span}\{e_1, \ldots, e_k\}$,
with $\mathbb{P}$-a.s.~continuous paths, such that, $\mathbb{P}$-a.s., the following identity holds
for all $t \in [0,T]$:
\[
  u^{(k)}(t) = \mathcal P_k u_0
  - \int_0^t \mathcal P_k\Bigl(
  \mathcal A u^{(k)}(s) + \bigl(u^{(k)}(s)\cdot\nabla\bigr)u^{(k)}(s) \Bigr)\,\mathrm{d}s
  + \sum_{n=1}^k \int_0^t \mathcal P_kL_{\zeta_n}u^{(k)}(s) \, \mathrm{d}W_n(s).
\]
Moreover, for any $v \in \operatorname{span}\{e_1, \ldots, e_k\}$,
the following equalities hold:
\begin{align*}
  &-2\Bigl\langle \mathcal P_k \Bigl( \mathcal A v + (v\cdot\nabla)v \Bigr), \, v \Bigr\rangle
  + \sum_{n=1}^k \|\mathcal P_kL_{\zeta_n}v\|_{\mathbb{L}_\sigma^2}^2 \\
  ={}&-2\Bigl\langle  \mathcal A v + (v\cdot\nabla)v, \, v \Bigr\rangle
  + \sum_{n=1}^k \|\mathcal P_kL_{\zeta_n}v\|_{\mathbb{L}_\sigma^2}^2 \\
  ={} & -2\|v\|_{X_1}^2 - \sum_{n=1}^k \|(I-\mathcal{P}_k)L_{\zeta_n}v\|_{\mathbb{L}_{\sigma}^2}^2 
  - \sum_{n=k+1}^\infty \|L_{\zeta_n}v\|_{\mathbb{L}_{\sigma}^2}^2, \\
   &\langle \mathcal P_k L_{\zeta_n}v, \, v\rangle
  = \langle (\zeta_n\cdot\nabla)v, \, v\rangle = 0, \qquad  1 \leqslant n \leqslant k,
\end{align*}
where $\langle \cdot, \cdot \rangle$ denotes the $\mathbb{L}^2$-inner product.
Applying It\^o's formula to $\|u^{(k)}(t)\|_{\mathbb{L}_\sigma^2}^2$ and using these relations, we obtain,
$\mathbb{P}$-a.s., for all $t\in[0,T]$,
\[
  \|u^{(k)}(t)\|_{\mathbb{L}_\sigma^2}^2
  +2\int_0^t\|u^{(k)}(s)\|_{X_1}^2\,\mathrm{d}s
  \leqslant \|\mathcal P_k u_0\|_{\mathbb{L}_\sigma^2}^2.
\]
Since $u_0\in L_{\mathcal F_0}^\infty(\Omega;\mathbb{L}_\sigma^2)$
and $\mathcal P_k$ is non-expansive on $\mathbb{L}_\sigma^2$, we deduce that
\begin{equation}
  \label{eq:un-basic}
  \sup_{k\geqslant1}\left(
    \|u^{(k)}\|_{L^\infty(\Omega;C([0,T];\mathbb{L}_\sigma^2))}
    +\|u^{(k)}\|_{L^\infty(\Omega;L^2(0,T;X_1))}
  \right)<\infty.
\end{equation}

\medskip\noindent
\textbf{Step 3.}
Fix $0 < \varrho < \varrho^{*}$.
In this step, we establish the following uniform bounds for the Galerkin
approximations $(u^{(k)})_{k\geqslant 1}$ under the additional assumption
$u_0 \in L_{\mathcal{F}_0}^2(\Omega;X_{1/2})$:
\begin{align}
\sup_{k \geqslant 1} \|u^{(k)}\|_{L^2(\Omega;C([0,T];X_{\varrho}))} &< \infty,
\label{eq:un-C} \\
\sup_{k \geqslant 1} \|u^{(k)}\|_{L^2(\Omega \times [0,T];X_{\varrho+1})} &< \infty.
\label{eq:un-L2}
\end{align}
The proof is organized into three substeps.

\textbf{(a).} 
For every $0 \leqslant \alpha \leqslant \varrho^{*}$, we establish that
the following inequality holds $\mathbb{P}$-a.s.~for all $t \in [0,T]$:
\begin{equation}
\label{eq:un-energy-ineq}
\begin{aligned}
& \|u^{(k)}(t)\|_{X_{\alpha}}^2 + c\int_0^t \|u^{(k)}(s)\|_{X_{\alpha+1}}^2  \mathrm{d}s \\
\leqslant{}& c\|\mathcal{P}_ku_0\|_{X_{\alpha}}^2
 + c \int_0^t \|u^{(k)}(s)\|_{X_{\alpha}}^2 \|u^{(k)}(s)\|_{X_{1}}^2  \mathrm{d}s \\
& + 2 \int_0^t \sum_{n=1}^k \bigl\langle  \mathcal{A}^{\alpha} u^{(k)}(s), \, \mathcal{P}_kL_{\zeta_n}u^{(k)}(s) \bigr\rangle \, \mathrm{d}W_n(s).
\end{aligned}
\end{equation}
To this end, applying It\^o's formula to $\|u^{(k)}(t)\|_{\mathcal{X}_{\alpha}}^2$
yields the following identity, which holds $\mathbb{P}$-a.s.~for all $t \in [0,T]$:
\begin{align*}
& \|u^{(k)}(t)\|_{\mathcal{X}_{\alpha}}^2
= \|\mathcal{P}_ku_0\|_{\mathcal{X}_{\alpha}}^2
 - 2\int_0^t \left( \|u^{(k)}(s)\|_{\mathcal{X}_{\alpha+1}}^2
 + \bigl\langle \mathcal{P}_k((u^{(k)}(s)\cdot\nabla)u^{(k)}(s)), \, \mathcal{A}^\alpha u^{(k)}(s) \bigr\rangle \right) \, \mathrm{d}s \\
& + \int_0^t \sum_{n=1}^{k} \|\mathcal{P}_kL_{\zeta_n}u^{(k)}(s)\|_{\mathcal{X}_{\alpha}}^2 \, \mathrm{d}s
 + 2 \int_0^t \sum_{n=1}^k \bigl\langle \mathcal{A}^{\alpha} u^{(k)}(s), \, \mathcal{P}_kL_{\zeta_n}u^{(k)}(s) \bigr\rangle \, \mathrm{d}W_n(s).
\end{align*}
The nonlinear term satisfies the following estimate:
\begin{align*}
2\big| \bigl\langle \mathcal{P}_k((u^{(k)}\cdot\nabla)u^{(k)}), \, \mathcal{A}^\alpha u^{(k)} \bigr\rangle \big|
&= 2\big| \bigl\langle \mathcal{P}((u^{(k)}\cdot\nabla)u^{(k)}), \, \mathcal{A}^\alpha u^{(k)} \bigr\rangle \big| \\
&\leqslant \epsilon \|u^{(k)}\|_{\mathcal{X}_{\alpha+1}}^2
+ \frac{c}{\epsilon} \|\mathcal{P}((u^{(k)}\cdot\nabla)u^{(k)})\|_{\mathcal{X}_{\alpha-1}}^2,
\quad \forall \epsilon > 0,
\end{align*}
where the equality holds because $\mathcal{A}^\alpha u^{(k)}$ lies in the range of
$\mathcal{P}_k$.
The It\^o correction term is bounded as follows:
\begin{align*}
\sum_{n=1}^k \|\mathcal{P}_kL_{\zeta_n}u^{(k)}\|_{\mathcal{X}_{\alpha}}^2
\leqslant \sum_{n=1}^\infty \|L_{\zeta_n}u^{(k)}\|_{\mathcal{X}_{\alpha}}^2 
= \|F(u^{(k)})\|_{\mathcal{L}_2(\ell^2, \mathcal{X}_{\alpha})}^2 
\leqslant C_F(\alpha) \|u^{(k)}\|_{\mathcal{X}_{\alpha+1}}^2,
\end{align*}
where $C_F(\alpha) := \|F\|_{\mathcal{L}(\mathcal{X}_{\alpha+1}, \, \mathcal{L}_2(\ell^2,\mathcal{X}_{\alpha}))}^2$;
here we used that $\|\mathcal{P}_k\|_{\mathcal{L}(\mathcal{X}_{\alpha}, \, \mathcal{X}_{\alpha})} =1$.
Combining these estimates with the fact that $\|\mathcal{P}_k\|_{\mathcal{L}(\mathcal{X}_{\alpha}, \mathcal{X}_{\alpha})}=1$,
we obtain
\begin{align*}
& \|u^{(k)}(t)\|_{\mathcal{X}_{\alpha}}^2
 + \int_0^t \big( 2 - \epsilon - C_F(\alpha) \big) \|u^{(k)}(s)\|_{\mathcal{X}_{\alpha+1}}^2 \, \mathrm{d}s \\
\leqslant{}& \|u_0\|_{\mathcal{X}_{\alpha}}^2
 + \frac{c}{\epsilon} \int_0^t \|\mathcal{P}((u^{(k)}(s)\cdot\nabla)u^{(k)}(s))\|_{\mathcal{X}_{\alpha-1}}^2 \, \mathrm{d}s \\
& + 2 \int_0^t \sum_{n=1}^k \bigl\langle \mathcal{A}^{\alpha} u^{(k)}(s), \, \mathcal{P}_kL_{\zeta_n}u^{(k)}(s) \bigr\rangle \, \mathrm{d}W_n(s).
\end{align*}
Since $C_F(\alpha) < 2$ (see \eqref{eq:PF}), we may choose $\epsilon > 0$ sufficiently small so that $2 - \epsilon - C_F(\alpha) > 0$.
Finally, combining the equivalence relation \eqref{eq:X-calX-equiv-negative} with the product estimate from \cite[Lemma 2.2]{Giga1985B},
\[
\|\mathcal{P}((u^{(k)}\cdot\nabla)u^{(k)})\|_{X_{\alpha-1}} \leqslant c \|u^{(k)}\|_{X_{\alpha}} \|u^{(k)}\|_{X_{1}},
\]
we obtain the desired inequality \eqref{eq:un-energy-ineq}.

\textbf{(b).} Let us prove \cref{eq:un-C}.
Apply the stochastic Gronwall lemma \cite[Theorem 4]{Scheutzow2013} to \eqref{eq:un-energy-ineq} with $\alpha = \varrho^{*}$. For any $p \in (0,1)$ and $\nu \in (1,1/p)$, this yields
\[
\|u^{(k)}\|_{L^{2p}(\Omega;C([0,T];X_{\varrho^{*}}))}
\leqslant c \left(\mathbb{E} \exp\left(\frac{p\nu}{\nu-1} \int_0^T \|u^{(k)}(t)\|_{X_{1}}^2  \mathrm{d}t \right) \right)^{\big(1-\frac{1}{\nu}\big)\frac{1}{2p}}
\|u_0\|_{L^{2p\nu}(\Omega;X_{\varrho^{*}})}.
\]
By the uniform bound \eqref{eq:un-basic}, the exponential moment on the right-hand side is finite and
uniformly bounded in $k$. Hence,
\[
\sup_{k \geqslant 1} \|u^{(k)}\|_{L^{2p}(\Omega;C([0,T];X_{\varrho^{*}}))}
\leqslant c \|u_0\|_{L^{2p\nu}(\Omega;X_{\varrho^{*}})},
\quad p \in (0,1), \, \nu \in (1,1/p).
\]
Since $u_0 \in L_{\mathcal{F}_0}^2(\Omega; X_{1/2}) \hookrightarrow L^{2p\nu}(\Omega;X_{\varrho^*})$,
we may take $p = \varrho / \varrho^{*} \in (0,1)$ to obtain
\[
\sup_{k \geqslant 1} \|u^{(k)}\|_{L^{2\varrho/\varrho^{*}}(\Omega;C([0,T];X_{\varrho^{*}}))} < \infty.
\]
Combined with \eqref{eq:un-basic}, this establishes the uniform bound \eqref{eq:un-C} via the interpolation inequality
\[
\|u^{(k)}\|_{L^2(\Omega;C([0,T];X_{\varrho}))}
\leqslant \|u^{(k)}\|_{L^\infty(\Omega;C([0,T];\mathbb{L}_{\sigma}^2))}^{1-\varrho/\varrho^{*}}
\|u^{(k)}\|_{L^{2\varrho/\varrho^{*}}(\Omega;C([0,T];X_{\varrho^{*}}))}^{\varrho/\varrho^{*}}.
\]

\textbf{(c).} Let us prove \cref{eq:un-L2}.
Taking the expectation of \eqref{eq:un-energy-ineq} with $\alpha = \varrho$
and $t = T$ and noting that the stochastic integral is a martingale with zero mean, we obtain
\begin{align*}
\mathbb{E} \Big[ \int_0^T \|u^{(k)}(t)\|_{X_{\varrho+1}}^2 \, \mathrm{d}t \Big]
&\leqslant c \mathbb{E} \left[
\|u_0\|_{X_{\varrho}}^2
+ \int_0^T \|u^{(k)}(t)\|_{X_{\varrho}}^2 \|u^{(k)}(t)\|_{X_{1}}^2 \, \mathrm{d}t 
\right] \\
& \leqslant c \mathbb{E}\left[
   \|u_0\|_{X_{\varrho}}^2 
   + \|u^{(k)}\|_{C([0,T];X_{\varrho})}^2
   \|u^{(k)}\|_{L^2(0,T;X_{1})}^2
   \right].
\end{align*}
The right-hand side is uniformly bounded in $k$ by \eqref{eq:un-basic},
\eqref{eq:un-C}, and the assumption $u_0 \in L_{\mathcal{F}_0}^2(\Omega;X_{1/2})$.
This establishes the uniform bound \eqref{eq:un-L2}.

\medskip\noindent
\textbf{Step 4.}
Let us prove \cref{eq:u-basic,eq:u-L2,eq:u-C}.
Fix $ 0 < \varrho < \varrho^*$ as in Step 3 and assume in addition that $u_0 \in L_{\mathcal{F}_0}^2(\Omega;X_{1/2})$.
In view of the uniform bounds \cref{eq:un-basic,eq:un-L2} and arguing as in \cite[Theorem~4.21]{Flandoli2023},
we conclude that the model problem \cref{eq:model} admits a unique solution $u$
in the sense of \cref{def:sol-def}, and \cref{eq:u-L2} holds in the sense that 
the $\mathbb{P}\otimes\mathrm{d}t$-equivalence class of $u$ belongs to $L^2(\Omega\times[0,T];X_{\varrho+1})$.
Moreover, \eqref{eq:u-basic} follows from the same argument as \eqref{eq:un-basic}.

It remains to establish \cref{eq:u-C}.
By the definition of the solution, it is standard
(see, e.g., \cite{Neerven2012b}) that, $\mathbb{P}$-a.s.,
the mild formulation holds for all $t \in [0,T]$:
\[
u(t) = e^{-t\mathcal{A}}u_0
- \int_0^t e^{(s-t)\mathcal{A}}\,
\mathcal{P}\bigl[(u(s)\cdot\nabla)u(s)\bigr] \,\mathrm{d}s
+ \int_0^t e^{(s-t)\mathcal{A}}\, F(u(s)) \,\mathrm{d}W(s),
\]
where $\bigl(e^{-t\mathcal{A}}\bigr)_{t \geqslant 0}$ is the analytic
semigroup on $\mathcal{X}_{-1}$ generated by $-\mathcal{A}$.
We show that each term on the right-hand side belongs to
$L^2(\Omega; C([0,T]; X_{\varrho}))$.

Since $u_0\in L_{\mathcal F_0}^2(\Omega;X_\varrho)$,
and since $(e^{-t\mathcal A})_{t\geqslant0}$ is a strongly continuous semigroup on
$\mathcal X_\varrho$, the equivalence relation \cref{eq:X-calX-equiv} with $\theta=\varrho$ implies
\[
  e^{-(\cdot)\mathcal A}u_0
  \in L^2\bigl(\Omega; C([0,T]; X_{\varrho})\bigr).
\]
A direct calculation using the Sobolev embedding theorem and H\"older's inequality yields
\[
  \bigl\|\mathcal{P}\bigl((v\cdot\nabla)v\bigr)\bigr\|_{X_{\varrho-1}}
  \leqslant c\, \|v\|_{\mathbb{L}_{\sigma}^2}\,\|v\|_{X_{\varrho+1}},
  \qquad v \in X_{\varrho+1}.
\]
Consequently, by H\"older's inequality and the established regularity properties \cref{eq:u-basic,eq:u-L2},
we obtain
\begin{equation}
  \label{eq:P-u-codt-nabla-u-regu}
  \mathcal{P}\bigl((u\cdot\nabla)u\bigr) \in L^2(\Omega \times [0,T]; X_{\varrho-1}).
\end{equation}
By the equivalence relation \cref{eq:X-calX-equiv-negative} (with $\theta=\varrho-1$),
and the standard $L^2$-maximal regularity estimate for analytic semigroups
(see, e.g., \cite[Theorem~3.3]{Neerven2012b}), we infer that
\[
  t \mapsto \int_0^t e^{(s-t)\mathcal{A}}\,
  \mathcal{P}\bigl((u\cdot\nabla)u\bigr) \,\mathrm{d}s
  \in L^2\bigl(\Omega; C([0,T]; \mathcal{X}_{\varrho})\bigr).
\]
Moreover, by \cref{eq:u-L2}, the equivalence relation \cref{eq:X-calX-equiv} (with $\theta=\varrho+1$),
and \cref{eq:PF} (with $\alpha=\varrho+1$), we obtain
\begin{equation}
  \label{eq:Fu-regu}
  F(u) \in L_{\mathbb{F}}^2(\Omega \times [0,T]; \mathcal{L}_2(\ell^2, \mathcal{X}_{\varrho})).
\end{equation}
Therefore, by \cite[Theorem~3.5]{Neerven2012b}, the process
\[
  t \mapsto \int_0^t e^{(s-t)\mathcal{A}}\, F(u(s)) \,\mathrm{d}W(s)
\]
has a continuous version in $L^2(\Omega; C([0,T]; \mathcal{X}_{\varrho}))$.
Combining the preceding conclusions and using the equivalence relation 
\cref{eq:X-calX-equiv} (with $\theta=\varrho$), we conclude \cref{eq:u-C}.
\end{proof}

\section{Spatial semidiscretization} \label{sec:spatial}
Let $\mathcal{K}_h$ be a conforming, quasi-uniform triangulation of the domain $\mathcal{O}$,
consisting of triangles. Denote by $h$ the maximum diameter of the elements in $\mathcal{K}_h$.
For each element $K \in \mathcal{K}_h$, define the local velocity space
\[
  \mathcal{V}(K) := \bigl[ P_1(K) \oplus \operatorname{span}\{\lambda_1 \lambda_2 \lambda_3\} \bigr]^2,
\]
where $P_1(K)$ denotes the space of linear polynomials on $K$,
and $\lambda_1, \lambda_2, \lambda_3$ are the barycentric (area) coordinates associated with the vertices of $K$.
The enrichment term $\lambda_1 \lambda_2 \lambda_3$ is the cubic bubble function vanishing on $\partial K$.
The classical MINI finite element spaces (cf.~\cite[Chapter~II, Section~4.1]{Girault1986}) are then defined as follows:
\begin{align*}
  \mathbb{L}_h &:= \left\{ v_h \in C(\overline{\mathcal{O}}; \mathbb{R}^2) : v_h|_K \in \mathcal{V}(K) \text{ for all } K \in \mathcal{K}_h,\; v_h|_{\partial\mathcal{O}} = 0 \right\}, \\
  M_h &:= \left\{ \phi_h \in C(\overline{\mathcal{O}}; \mathbb{R}) : \phi_h|_K \in P_1(K) \text{ for all } K \in \mathcal{K}_h \right\}.
\end{align*}
The discrete divergence-free subspace $\mathbb{L}_{\sigma,h} \subset \mathbb{L}_h$ is characterized by
the weak incompressibility condition
\begin{equation}
  \label{eq:L-sigma-h}
  \mathbb{L}_{\sigma,h} := \left\{
    v_h \in \mathbb{L}_h : \langle \nabla \cdot v_h, \, \phi_h \rangle = 0 \quad \forall \phi_h \in M_h
  \right\},
\end{equation}
where, here and throughout the sequel, $\langle\cdot, \cdot\rangle$ denotes
the $L^2(\mathcal{O})$-inner product for both scalar- and vector-valued
functions. The discrete Helmholtz projection $\mathcal{P}_h$
is defined as the $\mathbb{L}^2$-orthogonal projection onto $\mathbb{L}_{\sigma,h}$.
For every $\theta \in [0,1]$,
$\mathcal{P}_h$ admits an extension to $\mathcal{L}(\mathbb{H}_{D}^{-\theta},\, \mathbb{L}^{2})$,
still denoted by $\mathcal{P}_h$ for simplicity, which is defined via duality by
\[
  \langle \mathcal{P}_h v,\, v_h \rangle
  = \langle v,\, v_h \rangle_{\mathbb{H}_{D}^{-\theta},\, \mathbb{H}_{D}^{\theta}}
  \qquad \forall\, v \in \mathbb{H}_{D}^{-\theta},\ v_h \in \mathbb{L}_{\sigma,h}.
\]

The discrete Stokes operator $A_h : \mathbb{L}_{\sigma,h} \to \mathbb{L}_{\sigma,h}$ is given by
\[
\langle A_h u_h, v_h \rangle =
\int_{\mathcal{O}} \nabla u_h : \nabla v_h \, \mathrm{d}x, \quad \forall u_h, v_h \in \mathbb{L}_{\sigma,h}.
\]
The operator $A_h$ is self-adjoint and positive definite on $\mathbb{L}_{\sigma,h}$. Let $(X_{\theta,h})_{\theta\in\mathbb{R}}$ denote the Hilbert scale generated by $A_h$ on $\mathbb{L}_{\sigma,h}$; that is, for each $\theta\in\mathbb{R}$, $X_{\theta,h}$ is the space $\mathbb{L}_{\sigma,h}$ endowed with the norm
\[
\|v_h\|_{X_{\theta,h}} := \|A_h^{\theta/2}v_h\|_{\mathbb{L}^2},
\qquad \forall v_h\in\mathbb{L}_{\sigma,h}.
\]
In the particular case $\theta=1$, we have
\begin{equation}
  \label{eq:X1h}
  \|v_h\|_{X_{1,h}} = \|v_h\|_{\mathbb{H}_{D}^1}
  = \|\nabla v_h\|_{L^2(\mathcal{O};\mathbb{R}^{2\times2})},
  \quad \forall v_h \in X_{1,h}.
\end{equation}
The following estimate holds for all $v_h \in \mathbb{L}_{\sigma,h}$ and $\theta \in [0,1]$:
\begin{equation}
  \label{eq:Xh-to-HD}
  \|v_h\|_{\mathbb{H}_{D}^{\theta}} \leqslant \|v_h\|_{X_{\theta,h}}.
\end{equation}
Moreover, the inverse inequality
\begin{equation}
  \label{eq:inverse}
  \|v_h\|_{X_{\theta_2,h}} \leqslant c h^{\theta_1 - \theta_2} \|v_h\|_{X_{\theta_1,h}},
\end{equation}
holds for all $v_h \in \mathbb{L}_{\sigma,h}$ and all $-2 \leqslant \theta_1 \leqslant \theta_2 \leqslant 2$,
where $c > 0$ is independent of $h$ and $v_h$.

Assume that \cref{hyp:zeta} holds. For each integer $n\geqslant 1$,
define the operator $L_{\zeta_n,h}:\mathbb{L}^2 \to \mathbb{L}_{\sigma,h}$ by
\begin{equation}
  \label{eq:L-zeta-n-h}
  L_{\zeta_n,h}v := \mathcal{P}_h\bigl((\zeta_n\cdot\nabla)v\bigr),
  \qquad \forall v\in\mathbb{L}^2,
\end{equation}
where $(\zeta_n\cdot\nabla)v \in \mathbb{H}_{D}^{-1}$.
By integration by parts, the fact that $\mathcal{P}_h$ is the $\mathbb{L}^2$-orthogonal projection onto
$\mathbb{L}_{\sigma,h}$, the divergence-free property of $\zeta_n\in\mathbb{W}^{1,\infty}$
(cf.~\cref{hyp:zeta}), and the fact that both $v_h$ and $w_h$ vanish on $\partial\mathcal{O}$,
we have
\begin{equation}
  \label{eq:L-zeta-n-h-adjoint}
  \langle L_{\zeta_n,h}v_h, \, w_h\rangle
  =-\langle v_h, \, L_{\zeta_n,h}w_h\rangle
  \quad \text{for all } v_h,w_h \in \mathbb{L}_{\sigma,h}.
\end{equation}
Moreover, by the non-expansiveness of $\mathcal{P}_h$ on $\mathbb{L}^2$,
H\"older's inequality, and the identity
$\|v\|_{\mathbb{H}_{D}^1} = \|\nabla v\|_{L^2(\mathcal{O};\mathbb{R}^{2\times2})}$,
we obtain
\begin{equation}
  \label{eq:L-zeta-n-h-bound}
  \|L_{\zeta_n,h}v\|_{\mathbb{L}^2}
  \leqslant \|\zeta_n\|_{\mathbb{L}^\infty}
  \|v\|_{\mathbb{H}_D^1}
  \quad \text{for all } v\in\mathbb{H}_D^1.
\end{equation}

Combining \cref{eq:L-zeta-n-h-bound} with the identity
$\|\cdot\|_{X_{1,h}}=\|\cdot\|_{\mathbb H_D^1}$ on $X_{1,h}$,
we obtain
\[
  \sum_{n=1}^{\infty}\|L_{\zeta_n,h}v_h\|_{X_{0,h}}^{2}
  \leqslant
  \biggl(\sum_{n=1}^{\infty}
  \|\zeta_n\|_{\mathbb L^\infty}^{2}\biggr)
  \|v_h\|_{X_{1,h}}^{2}
  \qquad \text{for every } v_h \in X_{1,h}.
\]
Consequently, by the Cauchy--Schwarz inequality, the series
$\sum_{n=1}^{\infty}l_nL_{\zeta_n,h}v_h$ converges absolutely in
$X_{0,h}$ for every $v_h\in X_{1,h}$ and
$l=(l_n)_{n\geqslant 1}\in\ell^2$. The operator
$F_h \colon X_{1,h} \to \mathcal{L}_2(\ell^2, X_{0,h})$ is therefore
well defined by
\begin{equation}\label{eq:Fh-def}
  F_h(v_h)l
  := \sum_{n=1}^{\infty}l_nL_{\zeta_n,h}v_h,
  \qquad
  v_h\in X_{1,h},\;\; l=(l_n)_{n\geqslant 1}\in\ell^2,
\end{equation}
and, by the definition of the Hilbert--Schmidt norm,
\[
  \|F_h(v_h)\|_{\mathcal{L}_2(\ell^2,X_{0,h})}
  = \biggl(\sum_{n=1}^{\infty}
  \|L_{\zeta_n,h}v_h\|_{X_{0,h}}^{2}\biggr)^{1/2}
  \leqslant
  \biggl(\sum_{n=1}^{\infty}
  \|\zeta_n\|_{\mathbb L^\infty}^{2}\biggr)^{1/2}
  \|v_h\|_{X_{1,h}}
\]
holds for every $v_h \in X_{1,h}$. Consequently,
\begin{equation}\label{eq:Fh-bound0}
  \|F_h\|_{\mathcal{L}(X_{1,h},\mathcal{L}_2(\ell^2,X_{0,h}))}
  \leqslant
  \biggl(\sum_{n=1}^{\infty}
  \|\zeta_n\|_{\mathbb L^\infty}^{2}\biggr)^{1/2}.
\end{equation}

We consider the following spatial semidiscretization of the model problem \cref{eq:model}:
\begin{equation}
  \label{eq:uh}
  \begin{cases}
    \displaystyle
    \mathrm{d}u_h(t) = -\Big( A_h u_h(t) - \frac{1}{2}\sum_{n=1}^\infty L_{\zeta_n,h}^2 u_h(t)
    + \mathcal{P}_h G(u_h(t)) \Big) \, \mathrm{d}t \\
    \qquad\qquad\qquad {} + F_h(u_h(t)) \, \mathrm{d}W(t), \quad t \in [0,T], \\
    u_h(0) = \mathcal{P}_h u_0,
  \end{cases}
\end{equation}
where the nonlinear mapping $G$ is defined as
\begin{equation}
  \label{eq:G-def}
  G(v) := (v \cdot \nabla)v + \tfrac{1}{2} (\nabla \cdot v)v, \qquad \forall v \in \mathbb{H}_{D}^{1}.
\end{equation}
The key property of $G$ is that 
$\langle G(v), \, v \rangle = 0$ for all $v \in \mathbb{H}_{D}^1$.

We say that $u_h$ is a solution to the spatial semidiscretization \cref{eq:uh} if it satisfies the following conditions:
\begin{enumerate}[label=(\roman*)]
  \item[(i)] $u_h\colon \Omega\times[0,T]\to \mathbb{L}_{\sigma,h}$ is an $\mathbb{F}$-adapted process
  with $\mathbb{P}$-a.s.~continuous paths;
  \item[(ii)] $\mathbb{P}$-a.s., the following identity holds for all $t\in[0,T]$:
  \begin{equation}
    \label{eq:uh-def}
    \begin{aligned}
      u_h(t) &= \mathcal{P}_h u_0
      - \int_0^t \biggl( A_h u_h(s)
      - \frac{1}{2}\sum_{n=1}^\infty L_{\zeta_n,h}^2 u_h(s)
      + \mathcal{P}_h G(u_h(s)) \biggr) \, \mathrm{d}s \\
      &\quad {} + \int_0^t F_h(u_h(s)) \, \mathrm{d}W(s).
    \end{aligned}
  \end{equation}
\end{enumerate}
Two solutions are considered identical if they are indistinguishable.

\begin{remark}
  Observe that $\|\cdot\|_{\mathbb{H}_{D}^1} = \|\cdot\|_{X_{1,h}}$ on $X_{1,h}$.
  For each $v_h\in X_{0,h}$, two successive applications of \eqref{eq:L-zeta-n-h-bound},
  each followed by the inverse estimate \eqref{eq:inverse} with $(\theta_1,\theta_2)=(0,1)$, yield
  \[
    \|L_{\zeta_n,h}^2 v_h\|_{X_{0,h}}
    \leqslant c\, h^{-2}\, \|\zeta_n\|_{\mathbb{L}^\infty}^{2}
    \, \|v_h\|_{X_{0,h}},
  \]
  where $c>0$ is a constant independent of $h$, $n$, and $v_h$.
  Summing over $n$ and invoking \cref{eq:zeta-sum-cond}, we conclude that
  \[
    \sum_{n=1}^{\infty}L_{\zeta_n,h}^2v_h
    \quad \text{is absolutely convergent in $X_{0,h}$,}
  \]
  and
  \begin{equation}
    \label{eq:L-zeta-n-h-square}
    \Bigl\| \sum_{n=1}^{\infty}L_{\zeta_n,h}^2v_h \Bigr\|_{X_{0,h}}
    \leqslant ch^{-2} \Bigl(\sum_{n=1}^{\infty}\|\zeta_n\|_{\mathbb{L}^\infty}^2 \Bigr)\|v_h\|_{X_{0,h}}.
  \end{equation}
\end{remark}

We first present the well-posedness and stability of the spatial semidiscretization \cref{eq:uh}.

\begin{theorem}
  \label{thm:uh-stab}
  Assume that \cref{hyp:zeta} holds and that $u_0 \in L_{\mathcal{F}_0}^\infty(\Omega; \mathbb{L}_{\sigma}^2)$.
  Then the spatial semidiscretization \eqref{eq:uh} admits a unique solution $u_h$. Moreover,
  \begin{equation}
    \label{eq:uh-stab}
    \|u_h\|_{L^\infty(\Omega;  C([0,T]; \mathbb{L}^2))}
    + \|u_h\|_{L^\infty(\Omega; L^2(0,T; \mathbb{H}_{D}^1))}
    \leqslant (1+\sqrt{2}) \|u_0\|_{L^\infty(\Omega; \mathbb{L}^2)}.
  \end{equation}
\end{theorem}

\begin{proof}
  The proof follows from standard arguments; we provide a brief sketch for completeness.

\textbf{Existence and uniqueness.}
Since $X_{0,h}$ is finite-dimensional, the operator $A_h$ is bounded on it.
By \cref{eq:L-zeta-n-h-square} we have that $\sum_{n=1}^\infty L_{\zeta_n,h}^2$ is bounded on $X_{0,h}$,
and by \cref{eq:Fh-bound0} and the inverse estimate \cref{eq:inverse}
with $(\theta_1,\theta_2) = (0,1)$, we have that $F_h$ is a bounded linear operator from
$X_{0,h}$ to $\mathcal{L}_2(\ell^2,X_{0,h})$.
For any $v_h \in X_{0,h}$, the following identities hold:
\begin{align*}
  & \langle A_h v_h, v_h \rangle = \|v_h\|_{X_{1,h}}^2,
  \quad \langle \mathcal{P}_h G(v_h), v_h \rangle = 0, \\
  & \big\langle L_{\zeta_n,h}^2 v_h, v_h \big\rangle = -\|L_{\zeta_n,h} v_h\|_{X_{0,h}}^2
  \quad\text{for all $n\geqslant 1$,} \\
  & \|F_h(v_h)\|_{\mathcal{L}_2(\ell^2,X_{0,h})}^2 =
  \sum_{n=1}^{\infty} \|L_{\zeta_n,h}v_h\|_{X_{0,h}}^2.
\end{align*}
These yield the coercivity relation
\begin{equation}
  \label{eq:coercivity}
  2\Big\langle -A_h v_h + \tfrac{1}{2}\sum_{n=1}^\infty L_{\zeta_n,h}^2 v_h - \mathcal{P}_h G(v_h), \, v_h \Big\rangle
  + \|F_h(v_h)\|_{\mathcal{L}_2(\ell^2,X_{0,h})}^2
  + 2\|v_h\|_{X_{1,h}}^2 = 0.
\end{equation}
The nonlinear term $\mathcal{P}_h G$ is continuous on $X_{0,h}$. Moreover, for all $v_h, w_h \in \mathbb{L}_{\sigma,h}$, it satisfies the following growth and local monotonicity estimates:
\begin{align}
  \|\mathcal{P}_h G(v_h)\|_{X_{0,h}} &\leqslant c h^{-2} \|v_h\|_{X_{0,h}}^2, \label{eq:G-growth} \\
  \langle \mathcal{P}_h G(v_h) - \mathcal{P}_h G(w_h), \, v_h - w_h \rangle
  &\leqslant c h^{-2} \|w_h\|_{X_{0,h}} \|v_h - w_h\|_{X_{0,h}}^2, \label{eq:G-monotone}
\end{align}
where $c>0$ is a constant independent of $h$, $v_h$, and $w_h$.
These estimates follow from H\"older's inequality and standard inverse estimates
(see, e.g., \cite[Theorem~4.5.11]{Brenner2008}); the latter additionally relies on integration by parts.
The boundedness of the linear operators and the growth and local monotonicity conditions \eqref{eq:G-growth}--\eqref{eq:G-monotone} for the nonlinear term verify the hypotheses of \cite[Theorem~1.1]{LiuWei2010}. Consequently, \eqref{eq:uh} admits a unique solution $u_h$.

  \textbf{Uniform bound.}
The uniform bound for the solution is obtained via an argument analogous to \cite[Lemma~4.4]{Flandoli2023}.
First, we apply It\^o's formula to $\|u_h\|_{\mathbb{L}^2}^2$. The stochastic integral term vanishes due to the
skew-symmetry identity $\langle L_{\zeta_n,h}u_h, \, u_h \rangle = 0$ for all $n \geqslant 1$,
which follows from \cref{eq:L-zeta-n-h-adjoint}.
Hence, $\mathbb{P}$-a.s., the following identity holds for all $t \in [0,T]$:
\begin{align*}
  \|u_h(t)\|_{\mathbb{L}^2}^2 & = \|\mathcal{P}_hu_0\|_{\mathbb{L}^2}^2 
  + \int_0^t \Big\{ 2\Bigl\langle -A_hu_h(s) + \tfrac{1}{2} \sum_{n=1}^{\infty} L_{\zeta_n,h}^2 u_h(s) - \mathcal{P}_hG(u_h(s)), \, u_h(s) \Bigr\rangle \\
  & \qquad\qquad\qquad\qquad\qquad {} + \|F_h(u_h(s))\|_{\mathcal{L}_2(\ell^2,X_{0,h})}^2 \Big\} \, \mathrm{d}s .
\end{align*}
Next, substituting the coercivity relation \eqref{eq:coercivity} yields
that, $\mathbb{P}$-a.s., the following identity holds for all $t \in [0,T]$:
\[
  \|u_h(t)\|_{\mathbb{L}^2}^2 + 2 \int_0^t \|u_h(s)\|_{X_{1,h}}^2 \, \mathrm{d}s = \|\mathcal{P}_hu_0\|_{\mathbb{L}^2}^2.
\]
 Since $\mathcal{P}_h$ is non-expansive on $\mathbb{L}^2$,
 this immediately implies that, $\mathbb{P}$-a.s.,
  \[
    \|u_h\|_{C([0,T]; \mathbb{L}^2)}^2 + 2\|u_h\|_{L^2(0,T; X_{1,h})}^2 \leqslant 2\|\mathcal{P}_hu_0\|_{\mathbb{L}^2}^2
    \leqslant 2\|u_0\|_{\mathbb{L}^2}^2.
  \]
  Since $u_0 \in L_{\mathcal{F}_0}^{\infty}(\Omega;\mathbb{L}_{\sigma}^2)$
  and the identity $\|\cdot\|_{\mathbb{H}_{D}^1} = \|\cdot\|_{X_{1,h}}$ holds on $X_{1,h}$,
  the desired bound \cref{eq:uh-stab} follows.
  This completes the proof.
\end{proof}

Then, we present the convergence result for the spatial semidiscretization \cref{eq:uh}.

\begin{theorem}
  \label{thm:u-uh}
  Let \cref{hyp:zeta} hold, and let
  \[
    u_0 \in L_{\mathcal{F}_0}^{\infty}(\Omega; \mathbb{L}_{\sigma}^2)
    \cap L_{\mathcal{F}_0}^2(\Omega; X_{1/2}).
  \]
  Let $\varrho\in(0,\tfrac{1}{2})$ be the constant in \cref{prop:regu}.
  Let $u$ and $u_h$ denote the solutions of \eqref{eq:model} and of the
  spatial semidiscretization \cref{eq:uh}, respectively.
  Then there exists a constant $c>0$, depending on $\mathcal{O}$, $T$, $u_0$, $\varrho$,
  the quantity $\sum_{n=1}^\infty \|\zeta_n\|_{\mathbb{W}^{1,\infty}}^2$,
  the constant $\kappa_\zeta$ in \cref{eq:sigma-key},
  and the regularity parameters of the triangulation $\mathcal{K}_h$,
  but independent of the mesh size $h$, such that
  \begin{equation}
    \label{eq:u-uh}
    \|u - u_h\|_{L^2(\Omega; C([0,T]; \mathbb{L}^2))}
    + \|u - u_h\|_{L^2(\Omega\times(0,T); \mathbb{H}_{D}^1)}
    \leqslant c\, h^{\varrho} \log(1+1/h).
  \end{equation}
\end{theorem}

The remainder of this section is devoted to the proof of \cref{thm:u-uh}.
Throughout, $c>0$ denotes a generic constant that satisfies the same dependencies as in
\cref{thm:u-uh} but is, in addition, independent of the index $n \geqslant 1$ and 
the parameter $\epsilon > 0$; its value may vary from one occurrence to another.

\subsection{Preliminary estimates}
We begin by recalling the standard finite element approximation results and the associated operator-norm estimates.
The discrete Helmholtz projection $\mathcal{P}_h$ satisfies the classical approximation property
\begin{equation}
  \label{eq:Ph}
  \|I - \mathcal{P}_h\|_{\mathcal{L}(X_{\theta_2},\, \mathbb{H}_{D}^{\theta_1})}
  \leqslant c\, h^{\theta_2 - \theta_1},
  \qquad 0 \leqslant \theta_1 \leqslant 1, \ \ \theta_1 \leqslant \theta_2 \leqslant 2.
\end{equation}

By adapting the proof of \cite[Chapter~II, Theorem~4.1]{Girault1986},
we obtain the well-known approximation estimate
\begin{equation}
  \label{eq:Ah-approximation}
  \|A^{-1}\mathcal{P} - A_h^{-1}\mathcal{P}_h\|_{\mathcal{L}(\mathbb{H}_{D}^{-\theta}, \mathbb{L}^{2})}
  \leqslant c\, h^{2-\theta}, \qquad \theta \in [0,1].
\end{equation}
For $\theta \in [0,1]$, we estimate the composite operator $(I - \mathcal{P}_h)A^{-1}\mathcal{P}$ as follows:
\begin{align*}
  \|(I - \mathcal{P}_h)A^{-1}\mathcal{P}\|_{\mathcal{L}(\mathbb{H}_{D}^{-\theta}, \mathbb{L}^{2})}
  &\leqslant \|I - \mathcal{P}_h\|_{\mathcal{L}(X_{2-\theta}, \mathbb{L}^{2})}
  \|A^{-1}\|_{\mathcal{L}(X_{-\theta}, X_{2-\theta})}
  \|\mathcal{P}\|_{\mathcal{L}(\mathbb{H}_{D}^{-\theta}, X_{-\theta})} \\
  &\stackrel{\textup{(a)}}{\leqslant} c\, h^{2-\theta}
  \|\mathcal{P}\|_{\mathcal{L}(\mathbb{H}_{D}^{-\theta}, X_{-\theta})} \\
  &\stackrel{\textup{(b)}}{\leqslant} c\, h^{2-\theta},
\end{align*}
where step~(a) follows from \cref{eq:Ph} with
$(\theta_1,\theta_2) = (0,2-\theta)$, together with the fact that
$A^{-1}\colon X_{-\theta}\to X_{2-\theta}$ is an isometric isomorphism,
and step~(b) follows from \cref{lem:P}(iii) with $\alpha=-\theta$.
Combining this estimate with the approximation property
\cref{eq:Ah-approximation} and applying the triangle inequality yields
\begin{equation}
  \label{eq:PhA-AhPh}
  \|\mathcal{P}_h A^{-1}\mathcal{P} - A_h^{-1}\mathcal{P}_h\|
  _{\mathcal{L}(\mathbb{H}_{D}^{-\theta}, \mathbb{L}^{2})}
  \leqslant c\, h^{2-\theta},
  \qquad \theta \in [0,1].
\end{equation}
Setting $\theta=0$ in \cref{eq:PhA-AhPh}, restricting the
operators to $\mathbb{L}_{\sigma}^{2}$, and using the fact that
$\mathcal{P}$ acts as the identity on $\mathbb{L}_{\sigma}^2$,
we obtain
\begin{equation}
  \label{eq:PhA-AhPh-L2}
  \|\mathcal{P}_h A^{-1} - A_h^{-1}\mathcal{P}_h\|
  _{\mathcal{L}(\mathbb{L}_{\sigma}^{2}, \mathbb{L}^{2})}
  \leqslant c\, h^{2}.
\end{equation}

For each $v\in X_2$, it follows from \cref{eq:PhA-AhPh-L2} that
\[
  \|\mathcal{P}_h v-A_h^{-1}\mathcal{P}_hAv\|_{\mathbb{L}^2}
  =\|(\mathcal{P}_hA^{-1}-A_h^{-1}\mathcal{P}_h)Av\|_{\mathbb{L}^2}
  \leqslant c\,h^2\|Av\|_{\mathbb{L}_{\sigma}^2}
  = c\,h^2\|v\|_{X_2},
\]
where we have used the identity $\|Av\|_{\mathbb{L}_{\sigma}^2} = \|v\|_{X_2}$.
The inverse estimate \cref{eq:inverse} with $(\theta_1,\theta_2)=(0,2)$ then gives
\[
  \|\mathcal{P}_h v-A_h^{-1}\mathcal{P}_hAv\|_{X_{2,h}}
  \leqslant c\,\|v\|_{X_2}.
\]
Combining this with the identity
$\|A_h^{-1}\mathcal{P}_hAv\|_{X_{2,h}}=\|\mathcal{P}_hAv\|_{\mathbb{L}^{2}}$,
the $\mathbb{L}^{2}$-stability of $\mathcal{P}_h$, and the triangle inequality,
we obtain, for each $v \in X_2$,
\[
  \|\mathcal{P}_h v\|_{X_{2,h}}
  \leqslant \|\mathcal{P}_h v-A_h^{-1}\mathcal{P}_hAv\|_{X_{2,h}}
  +\|A_h^{-1}\mathcal{P}_hAv\|_{X_{2,h}}
  \leqslant c\,\|v\|_{X_2}.
\]
Moreover, the $\mathbb{L}^{2}$-stability estimate
$\|\mathcal{P}_h v\|_{X_{0,h}}\leqslant \|v\|_{X_0}$ holds for all $v\in X_0$.
Complex interpolation between the two preceding estimates yields
\begin{equation}
  \label{eq:Ph-interp}
  \|\mathcal{P}_h v\|_{X_{\theta,h}}
  \leqslant c\,\|v\|_{X_{\theta}}
  \qquad \text{for all } v\in X_{\theta} \text{ and } \theta\in[0,2].
\end{equation}

Combining the definition of the discrete divergence-free space
$\mathbb{L}_{\sigma,h}$ (see \cref{eq:L-sigma-h}) with the identity
\cref{eq:X1h}, we deduce that, for every $\alpha \in [0,2]$,
$v_h \in \mathbb{L}_{\sigma,h}$, and $\phi \in H^{\alpha}$,
\[
  |\langle \nabla\cdot v_h,\, \phi \rangle|
  = \min_{\phi_h \in M_h} |\langle \nabla\cdot v_h,\, \phi - \phi_h \rangle|
  \leqslant \sqrt{2} \|v_h\|_{X_{1,h}} \min_{\phi_h \in M_h} \|\phi - \phi_h\|_{L^2(\mathcal{O})}.
\]
Invoking the standard approximation properties of $M_h$, we then arrive at
\begin{equation}
  \label{eq:nabla-vh-phi}
  |\langle \nabla\cdot v_h,\, \phi \rangle|
  \leqslant c\, h^{\alpha}\, \|v_h\|_{X_{1,h}} \|\phi\|_{H^{\alpha}}
  \quad \forall\, v_h \in \mathbb{L}_{\sigma,h},\ \phi \in H^{\alpha}.
\end{equation}

The following lemma provides auxiliary bounds crucial for the subsequent analysis.

\begin{lemma}
\label{lem:auxi-ineqs}
The following bounds hold:
\begin{enumerate}
    \item[\textup{(i)}]
    For all $\theta \in [0,1]$,
    \begin{equation}
        \label{eq:Ph-I-P}
    \|\mathcal{P}_h(I - \mathcal{P})\|_{\mathcal{L}(\mathbb{H}^{\theta}, X_{-1,h})} \leqslant c h^{\theta + 1}.
    \end{equation}
    \item[\textup{(ii)}] For all $\theta \in (0,1/2)$,
    \begin{equation}
        \label{eq:P-Ph}
    \|\mathcal{P} - \mathcal{P}_h\|_{\mathcal{L}(\mathbb{H}^{\theta}, \mathbb{L}^2)} \leqslant c h^\theta.
    \end{equation}
\end{enumerate}
\end{lemma}

\begin{proof}
These estimates follow from standard finite element analysis for the Stokes problem; see, e.g., \cite[Chapter II]{Girault1986}.
For completeness, we present a self-contained proof.

  \textbf{Proof of (i).}
  Let $\theta \in [0,1]$ and $v \in \mathbb{H}^{\theta}$.
  By Theorem IV.3.5 in \cite{Boyer2013} and the boundedness of $ \mathcal{P} $ on $ \mathbb{H}^{\theta}  $ (cf.~\cref{lem:P}(i)),
  there exists a zero-mean function $\varphi \in H^{\theta+1}$ such that
  \begin{equation}
    \label{eq:varphi}
    v - \mathcal{P}v = \nabla\varphi, \qquad \|\varphi\|_{H^{\theta+1}} \leqslant c \|v\|_{\mathbb{H}^{\theta}}.
  \end{equation}
  Since $\mathcal{P}_h$ is the $\mathbb{L}^2$-orthogonal projection onto $\mathbb{L}_{\sigma,h}$, we have,
  for every $v_h \in \mathbb{L}_{\sigma,h}$,
  \[
    \langle v_h, \mathcal{P}_h(I - \mathcal{P})v \rangle
    = \langle v_h, (I - \mathcal{P})v \rangle
    = \langle v_h, \nabla \varphi \rangle
    = -\langle \nabla \cdot v_h, \varphi \rangle,
  \]
  where the last equality follows from integration by parts and 
  the fact that $v_h$ vanishes on $\partial\mathcal{O}$.
  Applying estimate \eqref{eq:nabla-vh-phi} with $\alpha = \theta+1$ and using the bound for $\|\varphi\|_{H^{\theta+1}}$ from \cref{eq:varphi} yields
  \[
    |\langle v_h, \mathcal{P}_h(I - \mathcal{P})v \rangle|
    \leqslant c h^{\theta+1} \|v_h\|_{X_{1,h}} \|\varphi\|_{H^{\theta+1}}
    \leqslant c h^{\theta+1} \|v_h\|_{X_{1,h}} \|v\|_{\mathbb{H}^{\theta}}.
  \]
  Taking the supremum over $v \in \mathbb{H}^{\theta}$ with $\|v\|_{\mathbb{H}^{\theta}}=1$ and $v_h \in \mathbb{L}_{\sigma,h}$ with $\|v_h\|_{X_{1,h}}=1$ gives
  the desired bound \cref{eq:Ph-I-P}.

\textbf{Proof of (ii).} Let $\theta \in (0,1/2)$. We decompose the difference as
\[
\mathcal{P} - \mathcal{P}_h = (I - \mathcal{P}_h)\mathcal{P} - \mathcal{P}_h(I - \mathcal{P})
\]
and estimate the two terms on the right-hand side.
For the first term, \cref{eq:Ph} with $\theta_1=0, \theta_2=\theta$ and the boundedness of $\mathcal{P}: \mathbb{H}^{\theta} \to X_{\theta}$ from Lemma~\ref{lem:P}(ii) imply
\[
\|(I - \mathcal{P}_h)\mathcal{P}\|_{\mathcal{L}(\mathbb{H}^{\theta}, \mathbb{L}^2)} \leqslant c h^{\theta}.
\]
For the second term, we employ the inverse inequality $\|w_h\|_{\mathbb{L}^2} \leqslant c h^{-1}\|w_h\|_{X_{-1,h}}$
for $w_h \in \mathbb{L}_{\sigma,h}$ (see~\eqref{eq:inverse}) together with Lemma~\ref{lem:auxi-ineqs}(i) to obtain
\[
\|\mathcal{P}_h(I - \mathcal{P})\|_{\mathcal{L}(\mathbb{H}^{\theta}, \mathbb{L}^2)} \leqslant c h^{-1} \|\mathcal{P}_h(I - \mathcal{P})\|_{\mathcal{L}(\mathbb{H}^{\theta},
X_{-1,h})} \leqslant c h^{\theta}.
\]
Combining these two bounds via the triangle inequality yields the desired estimate \cref{eq:P-Ph}.
\end{proof}

\subsection{The auxiliary operator \texorpdfstring{$\mathcal{J}_{h,\alpha}$}{Jh alpha}}

A principal difficulty in the numerical analysis of the model problem \eqref{eq:model}
stems from the limited regularity of the solution paths: in general,
they fail to belong to $X_{3/2+\varepsilon}$ for any $\varepsilon>0$ (see Remark~\ref{rem:regularity}).
Consequently, the expression $\mathcal{P}_h A u$ and the duality pairing
$\langle A u,\mathcal{P}v_h\rangle$ with $v_h\in\mathbb{L}_{\sigma,h}$ are not well-defined
in general. This regularity obstruction precludes the direct application of the techniques developed
in related works, including those in \cite{Breit2024}.

To circumvent this difficulty, we introduce the auxiliary operator
\begin{equation}
  \label{eq:Jh-def}
  \mathcal{J}_{h,\alpha} := A_h^{\alpha} \mathcal{P}_h A^{-\alpha},
  \qquad \alpha \in (0,1),
\end{equation}
which is given by the composition
\[
  X_{-2\alpha} \xrightarrow{\;A^{-\alpha}\;} \mathbb{L}_{\sigma}^{2}
  \xrightarrow{\;\mathcal{P}_h\;} X_{0,h}
  \xrightarrow{\;A_h^{\alpha}\;} X_{-2\alpha,h}.
\]
The operator $\mathcal{J}_{h,\alpha}$ has a regularizing effect: the fractional power $A^{-\alpha}$ first
maps an element of $X_{-2\alpha}$ into $\mathbb{L}_{\sigma}^{2}$;
the projection $\mathcal{P}_h$ then maps the result into the finite-dimensional space $X_{0,h}$;
and finally the discrete fractional power $A_h^{\alpha}$ is applied.
Thus $\mathcal{J}_{h,\alpha}$ serves as a discrete counterpart of the identity on $X_{-2\alpha}$.
Moreover, for every $\alpha \in (0,1)$ and
$\theta \in [-2\alpha,\, 2-2\alpha]$, invoking the fact that $A^{-\alpha}$
is an isometric isomorphism from $X_{\theta}$ onto $X_{\theta+2\alpha}$,
the stability estimate \cref{eq:Ph-interp} with $\theta$ replaced by
$\theta+2\alpha$, and the fact that $A_h^{\alpha}$ is an isometric
isomorphism from $X_{\theta+2\alpha,h}$ onto $X_{\theta,h}$, we conclude that
\begin{equation}
  \label{eq:Jh-stab}
  \|\mathcal{J}_{h,\alpha}\|_{\mathcal{L}(X_{\theta},X_{\theta,h})}
  \quad\text{is uniformly bounded with respect to $h$.}
\end{equation}
The key properties of $\mathcal{J}_{h,\alpha}$ are collected in the
following lemma.

\begin{lemma}
\label{lem:Jh}
The operator $ \mathcal{J}_{h,\alpha} $ satisfies the following estimates:
\begin{enumerate}
    \item[\textup{(i)}]  For all $\alpha \in (0,1)$, $ \theta_1 \in [-1,1] $, and $ \theta_2 \in [0,2] $,
    \begin{equation}
      \label{eq:Jh-Ph-0}
      \| \mathcal{J}_{h,\alpha} - \mathcal{P}_h \|_{\mathcal{L}(X_{\theta_2}, X_{\theta_1,h})}
      \leqslant ch^2 \int_0^\infty \frac{r^{-\alpha}}{1 + r^{\theta_2/2}}
      \min\!\left\{
        \frac{h^{(-2\alpha-\theta_1)\wedge 0}}{1+r^{0 \vee(-\alpha-\theta_1/2)}}, \, \frac{h^{-2-2\alpha-\theta_1}}{1 + r}
      \right\} \, \mathrm{d}r.
    \end{equation}

    \item[\textup{(ii)}]  For all $\alpha \in (0,1/2)$,
    \begin{align}
    & \|I - \mathcal{J}_{h,\alpha}\|_{\mathcal{L}(X_{\beta}, \mathbb{L}^2)} \leqslant c h^{\beta},
    \quad \beta \in [0,1], \label{eq:I-Jh-1} \\
    & \|I - \mathcal{J}_{h,\alpha}\|_{\mathcal{L}(X_{2-2\alpha}, \mathbb{H}_{D}^{1})} \leqslant c h^{1-2\alpha} \log(1+1/h). \label{eq:I-Jh-2} 
    \end{align}

    \item[\textup{(iii)}] For all $\alpha \in (0,1)$,
    \begin{equation}
      \label{eq:JhA-AhJh}
    \|\mathcal{J}_{h,\alpha} A - A_h \mathcal{J}_{h,\alpha}\|_{\mathcal{L}(X_{2-2\alpha},
    X_{-2,h})} \leqslant c h^{2-2\alpha}.
    \end{equation}

    \item[\textup{(iv)}] For all $\alpha \in (1/4,1/2)$,
    \begin{align}
    & \|\mathcal{J}_{h,\alpha}\mathcal{P} - \mathcal{P}_h\|_{\mathcal{L}(\mathbb{H}_{D}^{1-2\alpha}, X_{-1,h})} \leqslant c h^{2-2\alpha},
    \label{eq:JhP-Ph-1} \\
    & \|\mathcal{J}_{h,\alpha}\mathcal{P} - \mathcal{P}_h\|_{\mathcal{L}(\mathbb{H}_{D}^{-2\alpha}, X_{-1,h})}
     \leqslant c h^{1-2\alpha}. \label{eq:JhP-Ph-2} 
    \end{align}
\end{enumerate}
\end{lemma}

\begin{proof}
\textbf{(i)}  
Let $\alpha \in (0,1)$, $\theta_1 \in [-1,1]$, and $\theta_2 \in [0,2]$.  
From the definition of $\mathcal{J}_{h,\alpha}$ in~\cref{eq:Jh-def}, we obtain the identity
\[
  \mathcal{J}_{h,\alpha} - \mathcal{P}_h 
  = A_h^{\alpha} \bigl( \mathcal{P}_h A^{-\alpha} - A_h^{-\alpha} \mathcal{P}_h \bigr).
\]
Applying the Dunford--Taylor representation for fractional powers (cf.~Section 7.1 from Chapter 2 of \cite{Yagi2010}) yields
\[
  \mathcal{J}_{h,\alpha} - \mathcal{P}_h 
  = \frac{1}{2\pi \mathrm{i}} \int_{\Gamma} \lambda^{-\alpha} A_h^{\alpha}
  \Bigl[ \mathcal{P}_h (\lambda - A)^{-1} - (\lambda - A_h)^{-1} \mathcal{P}_h \Bigr] \, \mathrm{d}\lambda,
\]
where $\Gamma$ denotes the standard contour consisting of the rays $\{ r e^{\pm\mathrm{i}\pi/4} : r \in [0,\infty) \}$,
oriented with the positive real axis on its left.
This representation immediately implies the operator norm bound
\[
  \| \mathcal{J}_{h,\alpha} - \mathcal{P}_h \|_{\mathcal{L}(X_{\theta_2}, X_{\theta_1,h})}
  \leqslant \frac{1}{2\pi} \int_{\Gamma} |\lambda|^{-\alpha}
  \bigl\| A_h^{\alpha} \bigl[ \mathcal{P}_h (\lambda - A)^{-1} - (\lambda - A_h)^{-1} \mathcal{P}_h \bigr] \bigr\|_{\mathcal{L}(X_{\theta_2}, X_{\theta_1,h})} \, |\mathrm{d}\lambda|.
\]
To bound the integrand, we first observe the algebraic identity
\[
  A_h^{\alpha}\big[\mathcal{P}_h (\lambda - A)^{-1} - (\lambda - A_h)^{-1} \mathcal{P}_h \big]
  = (\lambda - A_h)^{-1} A_h^{1+\alpha} \bigl( A_h^{-1} \mathcal{P}_h - \mathcal{P}_h A^{-1} \bigr) A (\lambda - A)^{-1}.
\]
Using the approximation estimate
\[
  \|A_h^{-1}\mathcal{P}_h - \mathcal{P}_hA^{-1} \|_{\mathcal{L}(\mathbb{L}_{\sigma}^2, X_{0,h})} \leqslant c h^2
  \quad(\text{cf.~} \eqref{eq:PhA-AhPh-L2}),
\]
together with the resolvent bound
\[
  \| A (\lambda - A)^{-1} \|_{\mathcal{L}(X_{\theta_2}, \mathbb{L}_{\sigma}^2)} \leqslant \frac{c}{1 + |\lambda|^{\theta_2/2}},
\]
we obtain
\begin{align*}
  & \bigl\| A_h^{\alpha} \big[\mathcal{P}_h (\lambda - A)^{-1} - (\lambda - A_h)^{-1} \mathcal{P}_h\big] \bigr\|_{\mathcal{L}(X_{\theta_2},X_{\theta_1,h})} \\
  \leqslant{}& \frac{c h^2}{1 + |\lambda|^{\theta_2/2}} \,
  \|(\lambda - A_h)^{-1} A_h^{1+\alpha}\|_{\mathcal{L}(X_{0,h},X_{\theta_1,h})}.
\end{align*}
The operator $A_h$ is self-adjoint and positive definite on $\mathbb{L}_{\sigma,h}$.
Moreover, by \cref{eq:X1h} and the Poincar\'e inequality,
its smallest eigenvalue is bounded from below by a positive constant independent of $h$.
This implies the standard resolvent estimate
\[
  \| (\lambda - A_h)^{-1}A_h^{s} \|_{\mathcal{L}(X_{0,h}, X_{0,h})}
  \leqslant\frac{c}{1+|\lambda|^{1-s}}, \qquad s \in [0,1], \quad \lambda \in \Gamma.
\]
With the inverse inequality~\eqref{eq:inverse}, we further deduce
\begin{align*}
  \| (\lambda - A_h)^{-1} A_h^{1+\alpha} \|_{\mathcal{L}(X_{0,h}, X_{\theta_1,h})}
  &= \| (\lambda - A_h)^{-1} A_h^{1+\alpha+\theta_1/2} \|_{\mathcal{L}(X_{0,h}, X_{0,h})} \\
  &\leqslant c \min\!\left\{
    \frac{h^{(-2\alpha-\theta_1)\wedge 0}}{1+|\lambda|^{0 \vee(-\alpha-\theta_1/2)}}, \,
    \frac{h^{-2-2\alpha-\theta_1}}{1 + |\lambda|}
  \right\}.
\end{align*}
Consequently,
\begin{align*}
  & \bigl\| A_h^{\alpha} \big[\mathcal{P}_h (\lambda - A)^{-1} - (\lambda - A_h)^{-1} \mathcal{P}_h\big] \bigr\|_{\mathcal{L}(X_{\theta_2},X_{\theta_1,h})} \\
  &\qquad \leqslant \frac{c h^2}{1+|\lambda|^{\theta_2/2}} \min\!\left\{
     \frac{h^{(-2\alpha-\theta_1)\wedge 0}}{1+|\lambda|^{0 \vee(-\alpha-\theta_1/2)}}, \,
     \frac{h^{-2-2\alpha-\theta_1}}{1 + |\lambda|}
  \right\}.
\end{align*}
Substituting this bound into the Dunford--Taylor integral and parametrizing $\Gamma$ via $\lambda = r e^{\pm\mathrm{i}\pi/4}$ yields the desired estimate~\cref{eq:Jh-Ph-0}.

\textbf{(ii)}
Let $\alpha \in (0,\tfrac12)$ and $ \beta \in [0,1]$. Applying estimate \eqref{eq:Jh-Ph-0} with $(\theta_1,\theta_2) = (0,\beta)$ and $(1,2-2\alpha)$,
and employing the estimate \eqref{eq:Xh-to-HD} for $\theta \in \{0,1\}$, we obtain
\[
  \| \mathcal{J}_{h,\alpha} - \mathcal{P}_h \|_{\mathcal{L}(X_{\beta}, \mathbb{L}^2)} \leqslant c h^{\beta},
  \quad
  \| \mathcal{J}_{h,\alpha} - \mathcal{P}_h \|_{\mathcal{L}(X_{2-2\alpha}, \mathbb{H}_{D}^1)} \leqslant c h^{1-2\alpha} \log(1+1/h).
\]
Furthermore, the standard approximation properties of $\mathcal{P}_h$ (cf.~\eqref{eq:Ph}) provide
\[
  \|I - \mathcal{P}_h\|_{\mathcal{L}(X_{\beta},\mathbb{L}^2)} \leqslant c h^{\beta},
  \quad
  \|I - \mathcal{P}_h\|_{\mathcal{L}(X_{2-2\alpha},\mathbb{H}_{D}^{1})} \leqslant c h^{1-2\alpha}.
\]
Finally, the desired estimates \eqref{eq:I-Jh-1} and \eqref{eq:I-Jh-2} follow directly from the triangle inequality,
$\|I - \mathcal{J}_{h,\alpha}\| \leqslant \|I - \mathcal{P}_h\| + \|\mathcal{J}_{h,\alpha} - \mathcal{P}_h\|$,
in conjunction with the four bounds established above.

\textbf{(iii)}  
Fix $\alpha \in (0,1)$. From the definition \eqref{eq:Jh-def}, a direct computation yields the commutator identity
\[
\mathcal{J}_{h,\alpha} A - A_h \mathcal{J}_{h,\alpha}
= A_h^{\alpha+1} \bigl( A_h^{-1}\mathcal{P}_h - \mathcal{P}_h A^{-1} \bigr) A^{1-\alpha}.
\]
We estimate its operator norm from $X_{2-2\alpha}$ to $X_{-2,h}$ as follows:
\begin{align*}
\bigl\| \mathcal{J}_{h,\alpha} A - A_h \mathcal{J}_{h,\alpha} \bigr\|_{\mathcal{L}(X_{2-2\alpha},X_{-2,h})}
&= \bigl\| A_h^{\alpha+1}( A_h^{-1}\mathcal{P}_h - \mathcal{P}_h A^{-1} ) A^{1-\alpha} \bigr\|_{\mathcal{L}(X_{2-2\alpha},X_{-2,h})} \\
&= \bigl\| ( A_h^{-1}\mathcal{P}_h - \mathcal{P}_h A^{-1} ) A^{1-\alpha} \bigr\|_{\mathcal{L}(X_{2-2\alpha},X_{2\alpha,h})} \\
&\leqslant c h^{-2\alpha} \bigl\| ( A_h^{-1}\mathcal{P}_h - \mathcal{P}_h A^{-1} ) A^{1-\alpha} \bigr\|_{\mathcal{L}(X_{2-2\alpha},\mathbb{L}^2)} \\
&\leqslant c h^{-2\alpha} \bigl\|  A_h^{-1}\mathcal{P}_h - \mathcal{P}_h A^{-1} \bigr\|_{\mathcal{L}(\mathbb{L}_{\sigma}^2,\mathbb{L}^2)}
\bigl\| A^{1-\alpha} \bigr\|_{\mathcal{L}(X_{2-2\alpha},\mathbb{L}_{\sigma}^2)},
\end{align*}
where the first inequality uses the inverse inequality \eqref{eq:inverse}
with $(\theta_1,\theta_2) = (0,2\alpha)$.
Finally, using \cref{eq:PhA-AhPh-L2} and the fact that $A^{1-\alpha} \colon X_{2-2\alpha} \to \mathbb{L}_{\sigma}^2$
is an isometric isomorphism, we obtain \eqref{eq:JhA-AhJh}.

\textbf{(iv)}  
Let $\alpha \in (\tfrac14,\tfrac12)$. To prove \eqref{eq:JhP-Ph-1}, we write
\[
\mathcal{J}_{h,\alpha} \mathcal{P} - \mathcal{P}_h = (\mathcal{J}_{h,\alpha} - \mathcal{P}_h)\mathcal{P} + \mathcal{P}_h(\mathcal{P} - I).
\]
By Lemma~\ref{lem:P}(ii) and the continuous embedding 
$\mathbb{H}_{D}^{1-2\alpha} \hookrightarrow \mathbb{H}^{1-2\alpha}$,
the Helmholtz projection $\mathcal{P}$ is bounded from $\mathbb{H}_{D}^{1-2\alpha}$
to $X_{1-2\alpha}$, since $1-2\alpha \in (0,1/2)$.
Hence, applying estimate \eqref{eq:Jh-Ph-0} with $(\theta_1,\theta_2) = (-1,1-2\alpha)$ gives 
\[
\|(\mathcal{J}_{h,\alpha} - \mathcal{P}_h)\mathcal{P}\|_{\mathcal{L}(\mathbb{H}_{D}^{1-2\alpha},X_{-1,h})}
\leqslant c h^{2-2\alpha}.
\]
On the other hand, Lemma~\ref{lem:auxi-ineqs}(i) and the 
continuous embedding $\mathbb{H}_{D}^{1-2\alpha} \hookrightarrow \mathbb{H}^{1-2\alpha}$ imply
\[
\|\mathcal{P}_h(\mathcal{P} - I)\|_{\mathcal{L}(\mathbb{H}_{D}^{1-2\alpha},X_{-1,h})} \leqslant c h^{2-2\alpha}.
\]
Combining these two bounds yields \eqref{eq:JhP-Ph-1}.

To establish \eqref{eq:JhP-Ph-2}, we decompose
\begin{align*}
\mathcal{J}_{h,\alpha}\mathcal{P} - \mathcal{P}_h
&= A_h^\alpha \mathcal{P}_h A^{-\alpha}\mathcal{P} - \mathcal{P}_h \\
&= A_h^\alpha (\mathcal{P}_h - \mathcal{J}_{h,1-\alpha}) A^{-\alpha}\mathcal{P}
   + A_h \bigl( \mathcal{P}_h A^{-1}\mathcal{P} - A_h^{-1} \mathcal{P}_h \bigr).
\end{align*}
For the first term, observe that
\[
\|A_h^\alpha\|_{\mathcal{L}(X_{2\alpha-1,h}, X_{-1,h})} = 1,
\quad
\|A^{-\alpha}\|_{\mathcal{L}(X_{-2\alpha}, \mathbb{L}_{\sigma}^2)} = 1,
\]
and $\mathcal{P} \in \mathcal{L}(\mathbb{H}_{D}^{-2\alpha}, X_{-2\alpha})$ by Lemma~\ref{lem:P}(iii). Therefore,
\begin{align*}
& \|A_h^{\alpha}(\mathcal{P}_h - \mathcal{J}_{h,1-\alpha})A^{-\alpha}\mathcal{P}\|_{\mathcal{L}(\mathbb{H}_{D}^{-2\alpha}, X_{-1,h})} \\
\leqslant{} & \|A_h^{\alpha}\|_{\mathcal{L}(X_{2\alpha-1,h},X_{-1,h})}
\|\mathcal{P}_h - \mathcal{J}_{h,1-\alpha} \|_{\mathcal{L}(\mathbb{L}_{\sigma}^2,X_{2\alpha-1,h})}
\| A^{-\alpha}\|_{\mathcal{L}(X_{-2\alpha}, \mathbb{L}_{\sigma}^2)}
\| \mathcal{P}\|_{\mathcal{L}(\mathbb{H}_{D}^{-2\alpha}, X_{-2\alpha})} \\
\leqslant{} & c\,\|\mathcal{P}_h - \mathcal{J}_{h,1-\alpha}\|_{\mathcal{L}(\mathbb{L}_{\sigma}^2, X_{2\alpha-1,h})}.
\end{align*}
Applying \eqref{eq:Jh-Ph-0} with $\alpha$ replaced by $1-\alpha$, $\theta_1 = 2\alpha - 1$, and $\theta_2 = 0$
(noting that $2\alpha - 1 \in (-\tfrac12,0)$ for $\alpha \in (\tfrac14,\tfrac12)$) gives
\[
\|\mathcal{P}_h - \mathcal{J}_{h,1-\alpha}\|_{\mathcal{L}(\mathbb{L}_{\sigma}^2, X_{2\alpha-1,h})} \leqslant c h^{1-2\alpha}.
\]
Hence,  
\begin{equation}
\label{eq:tmp1}
\|A_h^\alpha(\mathcal{P}_h - \mathcal{J}_{h,1-\alpha})A^{-\alpha}\mathcal{P}\|_{\mathcal{L}(\mathbb{H}_{D}^{-2\alpha}, X_{-1,h})}
\leqslant c h^{1-2\alpha}.
\end{equation}
For the second term, note that
\[
\|A_h(\mathcal{P}_hA^{-1}\mathcal{P} - A_h^{-1}\mathcal{P}_h)\|_{\mathcal{L}(\mathbb{H}_{D}^{-2\alpha},X_{-1,h})}
= \|\mathcal{P}_hA^{-1}\mathcal{P} - A_h^{-1}\mathcal{P}_h\|_{\mathcal{L}(\mathbb{H}_{D}^{-2\alpha},X_{1,h})}.
\]
By the inverse inequality $\|v_h\|_{X_{1,h}} \leqslant c h^{-1} \|v_h\|_{\mathbb{L}^2}$
(cf.~\eqref{eq:inverse} with $\theta_1=0$ and $\theta_2=1$) and
the estimate $ \|\mathcal{P}_hA^{-1}\mathcal{P} - A_h^{-1}\mathcal{P}_h\|_{\mathcal{L}(\mathbb{H}_{D}^{-2\alpha},\mathbb{L}^2)}
\leqslant ch^{2-2\alpha} $ (cf.~\eqref{eq:PhA-AhPh} with $ \theta = 2\alpha$),
\[
\|\mathcal{P}_hA^{-1}\mathcal{P} - A_h^{-1}\mathcal{P}_h\|_{\mathcal{L}(\mathbb{H}_{D}^{-2\alpha},X_{1,h})}
\leqslant c h^{-1} \|\mathcal{P}_hA^{-1}\mathcal{P} - A_h^{-1}\mathcal{P}_h\|_{\mathcal{L}(\mathbb{H}_{D}^{-2\alpha},\mathbb{L}^2)}
\leqslant c h^{1-2\alpha}.
\]
Combining this bound with $\eqref{eq:tmp1}$ and applying the triangle inequality gives the desired bound $\eqref{eq:JhP-Ph-2}$.
\end{proof}

\subsection{Proof of Theorem~\ref{thm:u-uh}}
Since $A$ extends uniquely to an isometric isomorphism from $X_{\varrho+1}$ onto $X_{\varrho-1}$,
the regularity estimate \eqref{eq:u-L2} implies
\[
  Au\in L^{2}\bigl(\Omega\times[0,T]; X_{\varrho-1}\bigr).
\]
Similarly, applying \eqref{eq:L-zeta-n-theta-1-sum} with $\theta=\varrho$ and invoking \eqref{eq:u-L2} once again, we obtain
\[
  \sum_{n=1}^{\infty} L_{\zeta_n}^{2}u
  \in L^{2}\bigl(\Omega\times[0,T]; X_{\varrho-1}\bigr).
\]
In view of \eqref{eq:P-u-codt-nabla-u-regu} and the initial condition
\[
  u_0 \in L_{\mathcal{F}_0}^{\infty}(\Omega;\mathbb{L}_{\sigma}^{2})
  \hookrightarrow
  L_{\mathcal{F}_0}^{\infty}(\Omega;X_{\varrho-1}),
\]
it follows that the process
\[
  t \mapsto u_0-\int_0^t \Bigl(
    Au(s)-\frac{1}{2}\sum_{n=1}^{\infty}L_{\zeta_n}^{2}u(s)
    +\mathcal{P}\bigl((u(s)\cdot\nabla)u(s)\bigr)
  \Bigr)\,\mathrm{d}s
\]
belongs to $L^{2}\bigl(\Omega;C([0,T];X_{\varrho-1})\bigr)$.
Furthermore, by \eqref{eq:Fu-regu} and the continuous embedding
$\mathcal{X}_{\varrho}\hookrightarrow X_{\varrho}$
(see \eqref{eq:X-calX-equiv}), the stochastic integral
\[
  t\mapsto\int_0^t F(u(s))\,\mathrm{d}W(s)
\]
admits a version in $L^{2}(\Omega;C([0,T];X_{\varrho}))$.
Recalling the definition of $\mathcal{A}$ in \eqref{eq:calA-def} and
the continuous embedding $X_{\varrho} \hookrightarrow X_{\varrho-1}$, we infer that the right-hand side of \cref{eq:u-def} belongs to
$L^{2}(\Omega;C([0,T];X_{\varrho-1}))$, and, in view of \cref{eq:u-C} and the same embedding, so does the left-hand side.
On the other hand, since \cref{eq:u-def} holds in $\mathcal{X}_{-1}$ for every
$t\in[0,T]$, $\mathbb{P}$-a.s., and since $\mathcal{X}_{-1} = X_{-1}$ with equivalent norms
(see \eqref{eq:X-calX-equiv-negative} with $\theta=-1$),
it follows that, $\mathbb{P}$-a.s., the identity \cref{eq:u-def} holds in $X_{-1}$ for every $t\in[0,T]$.
The density of the canonical embedding $X_{1} \hookrightarrow X_{1-\varrho}$ implies
by duality that
two elements of $C([0,T];X_{\varrho-1})$ coincide in $C([0,T];X_{-1})$ also means they 
are identical in $C([0,T];X_{\varrho-1})$.
Consequently, we conclude that, $\mathbb{P}$-a.s., the identity
\begin{equation}\label{eq:u-strong}
  u(t)=u_0
  -\int_0^t\Bigl(
    Au(s)-\frac{1}{2}\sum_{n=1}^{\infty}L_{\zeta_n}^{2}u(s)
    +\mathcal{P}\bigl((u(s)\cdot\nabla)u(s)\bigr)
  \Bigr)\,\mathrm{d}s
  +\int_0^t F(u(s))\,\mathrm{d}W(s)
\end{equation}
holds in $X_{\varrho-1}$ for every $t\in[0,T]$.

We write $\mathcal{J}_h:=\mathcal{J}_{h,(1-\varrho)/2}$ for brevity.
Applying $\mathcal{J}_h$ to both sides of \cref{eq:u-strong} and setting
$ e_h := u_h - \mathcal{J}_h u $,
we conclude that, $\mathbb{P}$-a.s., 
\[
\begin{aligned}
  e_h(t) ={}& (\mathcal{P}_h - \mathcal{J}_h)u_0 - \int_0^t \Bigl(A_h e_h(s) + A_h \mathcal{J}_h u(s) - \mathcal{J}_h A u(s)\Bigr)\,\mathrm{d}s \\
  & + \frac{1}{2} \int_0^t \sum_{n=1}^\infty \big( L_{\zeta_n,h}^2 u_h(s) - \mathcal{J}_h L_{\zeta_n}^2 u(s) \big) \,\mathrm{d}s \\
  & + \int_0^t \Bigl[ \mathcal{J}_h \mathcal{P}((u(s)\cdot\nabla)u(s)) - \mathcal{P}_h G(u_h(s)) \Bigr]\,\mathrm{d}s \\
  & + \sum_{n=1}^\infty \int_0^t \Bigl[ L_{\zeta_n,h}u_h(s) - \mathcal{J}_h L_{\zeta_n}u(s) \Bigr]\,\mathrm{d}W_n(s),
\end{aligned}
\]
holds for every $t \in [0,T]$.
Applying It\^o's formula to $ \|e_h(t)\|_{\mathbb{L}^2}^2 $ yields,
$\mathbb{P}$-a.s., the following identity holds for all $t \in [0,T]$:
\begin{equation}
\label{eq:eh-sde}
\|e_h(t)\|_{\mathbb{L}^2}^2 = \|(\mathcal{P}_h-\mathcal{J}_h)u_0\|_{\mathbb{L}^2}^2
+ 2 \int_0^t \sum_{k=1}^4 I_k(s)  \,\mathrm{d}s + M(t),
\end{equation}
where the terms $I_k$ for $k=1,\dots,4$ and the process $M$ are defined by
\begin{align}
I_1 &:= \bigl\langle -A_h e_h - A_h \mathcal{J}_h u + \mathcal{J}_h A u, \, e_h \bigr\rangle, \label{eq:I1-def} \\
I_2 &:= \frac{1}{2} \sum_{n=1}^\infty \bigl\langle L_{\zeta_n,h}^2 u_h - \mathcal{J}_h L_{\zeta_n}^2 u,\, e_h \bigr\rangle, \label{eq:I2-def} \\
I_3 &:= \bigl\langle \mathcal{J}_h \mathcal{P}((u\cdot\nabla)u) - \mathcal{P}_h G(u_h),\, e_h \bigr\rangle, \label{eq:I3-def} \\
I_4 &:= \frac{1}{2} \sum_{n=1}^\infty \bigl\| L_{\zeta_n,h} u_h - \mathcal{J}_h L_{\zeta_n} u \bigr\|_{\mathbb{L}^2}^2, \label{eq:I4-def}
\end{align}
and
\begin{equation}
M(t) := 2  \sum_{n=1}^\infty \int_0^t \bigl\langle L_{\zeta_n,h}u_h(s)
 - \mathcal{J}_h L_{\zeta_n}u(s),\, e_h(s) \bigr\rangle \,\mathrm{d}W_n(s). \label{eq:M-def}
\end{equation}


The proof proceeds in six steps. In Steps 1-5, we derive bounds for the terms $I_1, \dots, I_4$
and the process $M(t)$.
Finally, in Step 6, we combine these estimates with Gronwall's inequality
to establish the convergence rate stated in \cref{eq:u-uh}.

\medskip\noindent
\textbf{Step 1: Estimate of $I_1$.} We establish the bound
\begin{equation}
\label{eq:I1-bound}
I_1 \leqslant -\frac{3}{4} \|e_h\|_{\mathbb{H}_{D}^1}^2 + c h^{2\varrho} \|u\|_{X_{\varrho+1}}^2.
\end{equation}
By the definition of $I_1$ in \cref{eq:I1-def}, we decompose $I_1$ as
\[
I_1 = -\|e_h\|_{X_{1,h}}^2 + \langle \mathcal{J}_h A u - A_h \mathcal{J}_h u, \, e_h \rangle.
\]
To bound the second term, we argue as follows:
\begin{align*}
\langle \mathcal{J}_h A u - A_h \mathcal{J}_h u, \, e_h \rangle 
&\leqslant \|e_h\|_{X_{2,h}} \|\mathcal{J}_h A u - A_h \mathcal{J}_h u\|_{X_{-2,h}}  \\
&\leqslant c h^{-1} \|e_h\|_{X_{1,h}} \|\mathcal{J}_h A u - A_h \mathcal{J}_h u\|_{X_{-2,h}} \\
&\leqslant c h^{\varrho} \|e_h\|_{X_{1,h}} \|u\|_{X_{\varrho+1}},
\end{align*}
where the second inequality follows from the inverse inequality \cref{eq:inverse} with $(\theta_1,\theta_2) = (1,2)$, and the third from \cref{eq:JhA-AhJh} with $\alpha = (1-\varrho)/2$.
Substituting this estimate into the decomposition of $I_1$ and using the identity
$\|e_h\|_{X_{1,h}} = \|e_h\|_{\mathbb{H}_{D}^1}$ yields
\[
I_1 \leqslant -\|e_h\|_{\mathbb{H}_{D}^1}^2 + c h^\varrho \|e_h\|_{\mathbb{H}_{D}^1} \|u\|_{X_{\varrho+1}}.
\]
Applying Young's inequality then gives the desired bound \cref{eq:I1-bound}.

\medskip\noindent
\textbf{Step 2. (Estimate of $I_2$).} We establish the following bound for $I_2$:
\begin{equation}
\label{eq:I2-bound}
I_2 \leqslant
\frac{3}{8} \|e_h\|_{\mathbb{H}_{D}^1}^2
- \frac{1}{2} \sum_{n=1}^\infty \|L_{\zeta_n,h}e_h\|_{\mathbb{L}^2}^2
+ c h^{2\varrho} (\log(1+1/h))^2 \|u\|_{X_{\varrho+1}}^2.
\end{equation}
From the definition \cref{eq:I2-def}, we decompose $I_2$ as
\[
I_2 = \frac{1}{2}\sum_{n=1}^{\infty}\left\langle (L_{\zeta_n,h} - \mathcal{J}_hL_{\zeta_n})L_{\zeta_n}u, \, e_h \right\rangle 
 + \frac{1}{2}\sum_{n=1}^{\infty}\left\langle L_{\zeta_n,h} \big( L_{\zeta_n,h}u_h - L_{\zeta_n}u\big), \, e_h \right\rangle.
\]
Using the definition \cref{eq:L-zeta-n-h}, 
the fact that $\mathcal{P}_h$ is the $\mathbb{L}^2$-projection onto $\mathbb{L}_{\sigma,h}$,
integration by parts, the fact that each $\zeta_n \in \mathbb{W}^{1,\infty}$ is divergence-free
(cf.~\cref{hyp:zeta}), and the fact that $e_h$ vanishes on $\partial\mathcal{O}$,
the second term can be rewritten as 
\begin{align*}
  & \frac{1}{2}\sum_{n=1}^{\infty}\left\langle L_{\zeta_n,h} \big( L_{\zeta_n,h}u_h - L_{\zeta_n}u\big), \, e_h \right\rangle \\
   ={} &
  -\frac{1}{2}\sum_{n=1}^{\infty}\left\langle L_{\zeta_n,h}u_h - L_{\zeta_n}u, \, (\zeta_n \cdot \nabla)e_h \right\rangle \\
   ={} &
  \frac{1}{2}\sum_{n=1}^{\infty}\left\langle
     (L_{\zeta_n} - L_{\zeta_n,h})u
     + L_{\zeta_n,h}(I - \mathcal{J}_h)u 
     - L_{\zeta_n,h}e_h, \, (\zeta_n \cdot \nabla)e_h
     \right\rangle \\
    ={} &
    \frac{1}{2} \sum_{n=1}^\infty \langle (L_{\zeta_n} - L_{\zeta_n,h})u, \, (\zeta_n\cdot\nabla)e_h \rangle
    + \frac{1}{2} \sum_{n=1}^\infty \langle L_{\zeta_n,h}(I - \mathcal{J}_h)u, \, L_{\zeta_n,h}e_h \rangle \\
    & \qquad {} - \frac{1}{2} \sum_{n=1}^\infty \|L_{\zeta_n,h}e_h\|_{\mathbb{L}^2}^2.
\end{align*}
Hence,
\begin{equation}
\label{eq:I2-decomp}
I_2 = I_2^{(1)} + I_2^{(2)} + I_2^{(3)} - \frac{1}{2}\sum_{n=1}^\infty \|L_{\zeta_n,h}e_h\|_{\mathbb{L}^2}^2,
\end{equation}
where
\begin{align*}
I_2^{(1)} &:= \frac{1}{2}\sum_{n=1}^{\infty}\left\langle (L_{\zeta_n,h} - \mathcal{J}_hL_{\zeta_n})L_{\zeta_n}u, \, e_h \right\rangle, \\
I_2^{(2)} &:= \frac{1}{2}\sum_{n=1}^{\infty}\left\langle (L_{\zeta_n} - L_{\zeta_n,h})u, \, (\zeta_n \cdot \nabla) e_h \right\rangle, \\
I_2^{(3)} &:= \frac{1}{2}\sum_{n=1}^{\infty}\left\langle L_{\zeta_n,h}(I-\mathcal{J}_h)u, \, L_{\zeta_n,h} e_h \right\rangle.
\end{align*}
The desired bound \cref{eq:I2-bound} is then obtained by combining the estimates
\cref{eq:I2-1,eq:I2-2,eq:I2-3} for $ I_2^{(1)} $, $ I_2^{(2)} $, and $ I_2^{(3)} $,
 which are derived in parts (a)–(c) below.

\textit{Part (a): Estimate for \(I_2^{(1)}\).}
We estimate \(I_2^{(1)}\) as follows:
\[
  \begin{aligned}
    I_2^{(1)} &\stackrel{\text{(i)}}{=} \frac{1}{2} \sum_{n=1}^{\infty} \left\langle
    (\mathcal{P}_h - \mathcal{J}_h\mathcal{P})[(\zeta_n \cdot \nabla)L_{\zeta_n}u], \, e_h
    \right\rangle \\
    &\stackrel{\text{(ii)}}{\leqslant} \frac{1}{2} \sum_{n=1}^{\infty} \left\| (\mathcal{P}_h - \mathcal{J}_h\mathcal{P})[(\zeta_n \cdot \nabla)L_{\zeta_n}u] \right\|_{X_{-1,h}} \left\| e_h \right\|_{X_{1,h}} \\
    &\stackrel{\text{(iii)}}{\leqslant} c h^\varrho \sum_{n=1}^{\infty} \left\| (\zeta_n \cdot \nabla)L_{\zeta_n}u \right\|_{\mathbb{H}_{D}^{\varrho-1}} \left\| e_h \right\|_{X_{1,h}}  \\
    &\stackrel{\text{(iv)}}{\leqslant} c h^\varrho \sum_{n=1}^{\infty} \| \zeta_n \|_{\mathbb{W}^{1,\infty}}^2 \left\| u \right\|_{X_{\varrho+1}} \left\| e_h \right\|_{X_{1,h}},
  \end{aligned}
\]
where (i) follows from the definitions in \cref{eq:L-zeta-n,eq:L-zeta-n-h}; (ii) follows from the Cauchy--Schwarz inequality; (iii) follows from \cref{eq:JhP-Ph-2} with \(\alpha=(1-\varrho)/2\); and (iv) follows from \cref{eq:zetan-cdot-nabla-v,eq:L-zeta-n-theta} with \(\theta=\varrho\).
Using the summability condition \cref{eq:zeta-sum-cond} together with the identity \(\|e_h\|_{X_{1,h}} = \|e_h\|_{\mathbb{H}_{D}^1}\), we obtain
\[
  I_2^{(1)} \leqslant c h^\varrho
  \left\| u \right\|_{X_{\varrho+1}} \left\| e_h \right\|_{\mathbb{H}_{D}^1}.
\]
Applying Young's inequality then yields
\begin{equation}
  \label{eq:I2-1}
  I_2^{(1)} \leqslant \frac{1}{8} \|e_h\|_{\mathbb{H}_{D}^1}^2
  + c h^{2\varrho} \|u\|_{X_{\varrho+1}}^2.
\end{equation}

\textit{Part (b): Estimate of $ I_2^{(2)} $.}
By the definitions \cref{eq:L-zeta-n,eq:L-zeta-n-h}, Cauchy--Schwarz inequality,
and \cref{lem:auxi-ineqs}(ii), we deduce
\begin{align*}
I_2^{(2)}
&= \frac{1}{2} \sum_{n=1}^{\infty} \langle (\mathcal{P} - \mathcal{P}_h)((\zeta_n\cdot\nabla)u), \, (\zeta_n \cdot \nabla)e_h \rangle \\
& \leqslant c h^{\varrho} \sum_{n=1}^{\infty}
\| (\zeta_n \cdot \nabla)u \|_{\mathbb{H}^{\varrho}} \| (\zeta_n \cdot \nabla)e_h \|_{\mathbb{L}^2}.
\end{align*}
Applying the standard product estimates
\begin{align*}
& \| (\zeta_n \cdot \nabla)u \|_{\mathbb{H}^{\varrho}}
\leqslant c \| \zeta_n \|_{\mathbb{W}^{1,\infty}} \| u \|_{X_{\varrho+1}}, \\
& \| (\zeta_n \cdot \nabla)e_h \|_{\mathbb{L}^2}
\leqslant c \| \zeta_n \|_{\mathbb{L}^{\infty}} \| e_h \|_{\mathbb{H}_{D}^1},
\end{align*}
then yields
\[
I_2^{(2)} \leqslant c h^{\varrho} \sum_{n=1}^{\infty} \| \zeta_n \|_{\mathbb{W}^{1,\infty}}^2
\| u \|_{X_{\varrho+1}} \| e_h \|_{\mathbb{H}_{D}^1}.
\]  
By the summability condition \eqref{eq:zeta-sum-cond} and applying Young's inequality, we arrive at
\begin{equation}  
\label{eq:I2-2}  
I_2^{(2)} \leqslant \frac{1}{8} \| e_h \|_{\mathbb{H}_{D}^1}^2
+ c h^{2\varrho} \| u \|_{X_{\varrho+1}}^2.  
\end{equation}

\textit{Part (c): Estimate of $ I_2^{(3)} $.}
We bound $I_2^{(3)}$ as follows:
\begin{align*}
  I_2^{(3)}
  &\leqslant \frac{1}{2} \sum_{n=1}^{\infty}
  \|L_{\zeta_n,h}(I-\mathcal{J}_h)u\|_{\mathbb{L}^2}
  \, \|L_{\zeta_n,h}e_h\|_{\mathbb{L}^2} \\
  &\leqslant c \sum_{n=1}^{\infty} \|\zeta_n\|_{\mathbb{L}^{\infty}}^{2}
  \, \|(I - \mathcal{J}_h)u\|_{\mathbb{H}_{D}^{1}} \, \|e_h\|_{\mathbb{H}_{D}^{1}} \\
  &\leqslant c\, h^{\varrho} \log(1 + 1/h) \left(\sum_{n=1}^{\infty} \|\zeta_n\|_{\mathbb{L}^{\infty}}^{2}\right)
  \|u\|_{X_{\varrho+1}} \, \|e_h\|_{\mathbb{H}_{D}^{1}} \\
  &\leqslant c\, h^{\varrho} \log(1 + 1/h) \|u\|_{X_{\varrho+1}} \, \|e_h\|_{\mathbb{H}_{D}^{1}},
\end{align*}
where the first inequality follows from the Cauchy--Schwarz inequality,
the second from \cref{eq:L-zeta-n-h-bound},
the third from \cref{eq:I-Jh-2} with $\alpha = (1-\varrho)/2$,
and the last from the summability condition \eqref{eq:zeta-sum-cond}.
An application of Young's inequality then yields
\begin{equation}
  \label{eq:I2-3}
  I_2^{(3)} \leqslant
  \frac{1}{8} \|e_h\|_{\mathbb{H}_{D}^{1}}^{2}
  + c\, h^{2\varrho} \bigl(\log(1 + 1/h)\bigr)^{2} \|u\|_{X_{\varrho+1}}^{2}.
\end{equation}

\medskip\noindent
\textbf{Step 3. (Estimate of $I_3$).} We establish the following bound for $I_3$:
\begin{equation}
\label{eq:I3-bound}
\begin{aligned}
  I_3 &\leqslant
  \frac{1}{8}\|e_h\|_{\mathbb{H}_{D}^1}^2
  + c \Big( \|u\|_{\mathbb{L}^2}^2\|u\|_{X_{1}}^2  +
  \|u_h\|_{\mathbb{L}^2}^2\|u_h\|_{\mathbb{H}_{D}^1}^2 \Big) \|e_h\|_{\mathbb{L}^2}^2 \\
  & \quad {} + ch^{2\varrho} (1+\|u\|_{\mathbb{L}^2}^2)\|u\|_{X_{\varrho}}^2 \|u\|_{X_{1}}^2
  + ch^{2\varrho} (\log(1+1/h))^2 \|u\|_{X_{\varrho+1}}^2.
\end{aligned}
\end{equation}
From the definition of $I_3$ in \cref{eq:I3-def}, the definition of $G$ in \cref{eq:G-def},
and the incompressibility of $ u $, we decompose $ I_3 $ as
\begin{align*}
I_3 & = \underbrace{\langle(\mathcal{J}_h\mathcal{P} - \mathcal{P}_h)((u\cdot\nabla)u), \, e_h\rangle}_{I_3^{(1)}} \\
& \qquad {} + \underbrace{\langle((u - u_h)\cdot\nabla) u, \, e_h\rangle
 + \langle(u_h\cdot\nabla)(u - u_h), \, e_h\rangle
 + \frac{1}{2}\langle\nabla\cdot(u - u_h), \, u_h \cdot e_h\rangle }_{I_3^{(2)}},
\end{align*}
where we used that $\mathcal{P}_h$ is the $\mathbb{L}^2$-orthogonal projection onto $\mathbb{L}_{\sigma,h}$.
For $ I_3^{(1)} $, using the Cauchy--Schwarz inequality and \cref{eq:JhP-Ph-2} with $ \alpha = (1-\varrho)/2 $, we have
\begin{align*}
I_3^{(1)} &\leqslant \|(\mathcal{J}_h\mathcal{P} - \mathcal{P}_h)((u\cdot\nabla)u)\|_{X_{-1,h}} \|e_h\|_{X_{1,h}} \\
&\leqslant c h^{\varrho} \|(u\cdot\nabla)u\|_{\mathbb{H}_{D}^{\varrho-1}} \|e_h\|_{X_{1,h}}.
\end{align*}
A duality argument, employing H\"older's inequality,
the continuous embeddings \cref{eq:X-to-HD} for $\theta \in \{\varrho,1\}$,
and the Sobolev embedding theorem, gives
\[
\|(u\cdot\nabla)u\|_{\mathbb{H}_{D}^{\varrho-1}} \leqslant c \|u\|_{X_{\varrho}} \|u\|_{X_{1}},
\]
and hence,
\[
I_3^{(1)} \leqslant c h^{\varrho} \|u\|_{X_{\varrho}} \|u\|_{X_{1}} \|e_h\|_{X_{1,h}}.
\]
Applying the equality $\|e_h\|_{X_{1,h}} = \|e_h\|_{\mathbb{H}_{D}^1}$ and Young's inequality then gives
\[
I_3^{(1)} \leqslant \frac{1}{16} \|e_h\|_{\mathbb{H}_{D}^1}^2 + c h^{2\varrho} \|u\|_{X_{\varrho}}^2 \|u\|_{X_{1}}^2.
\]
We now turn to $ I_3^{(2)} $. Integration by parts provides
\begin{align*}
I_3^{(2)} &= \langle ((u_h - u)\cdot\nabla) e_h, \, u\rangle + \langle \nabla\cdot(u_h - u), \, u \cdot e_h\rangle \\
&\quad + \langle(u_h\cdot\nabla)(u - u_h), \, e_h\rangle + \frac{1}{2}\langle\nabla\cdot(u - u_h), \, u_h \cdot e_h\rangle.
\end{align*}
By H\"older's inequality, it follows that
\begin{align*}
I_3^{(2)} \leqslant c\|u_h-u\|_{\mathbb{L}^4} \|e_h\|_{\mathbb{H}_{D}^1} \|u\|_{\mathbb{L}^4}
+ c (\|u\|_{\mathbb{L}^4} + \|u_h\|_{\mathbb{L}^4}) \|u-u_h\|_{\mathbb{H}_{D}^1} \|e_h\|_{\mathbb{L}^4}.
\end{align*}
Using Young's inequality, the identity $ u_h - u = e_h - (u - \mathcal{J}_h u) $,
the Gagliardo--Nirenberg interpolation inequality
$ \|v\|_{\mathbb{L}^4} \leqslant c\|v\|_{\mathbb{L}^2}^{1/2} \|v\|_{\mathbb{H}_{D}^1}^{1/2} $
for any $ v \in \mathbb{H}_{D}^{1} $ (see, e.g., \cite[Proposition III.2.35]{Boyer2013}),
and the continuous embedding $X_1 \hookrightarrow \mathbb{H}_{D}^1$ (see \cref{eq:X-to-HD} with $\theta=1$),
we deduce
\[
\begin{aligned}
I_3^{(2)} &\leqslant \frac{1}{16} \|e_h\|_{\mathbb{H}_{D}^1}^2
+ c \Big( \|u\|_{\mathbb{L}^2}^2\|u\|_{X_{1}}^2 + \|u_h\|_{\mathbb{L}^2}^2\|u_h\|_{\mathbb{H}_{D}^1}^2 \Big) \|e_h\|_{\mathbb{L}^2}^2 \\
&\quad + c \|u\|_{\mathbb{L}^2}^2\|u\|_{X_{1}}^2 \|u - \mathcal{J}_h u\|_{\mathbb{L}^2}^2 + c \|u - \mathcal{J}_h u\|_{\mathbb{H}_{D}^{1}}^2.
\end{aligned}
\]
Applying \cref{eq:I-Jh-1,eq:I-Jh-2} with $ \alpha = (1-\varrho)/2 $ and $ \beta = 1-2\alpha $ then yields
\[
\begin{aligned}
I_3^{(2)} &\leqslant \frac{1}{16} \|e_h\|_{\mathbb{H}_{D}^1}^2 + c \Big( \|u\|_{\mathbb{L}^2}^2\|u\|_{X_{1}}^2 + \|u_h\|_{\mathbb{L}^2}^2\|u_h\|_{\mathbb{H}_{D}^1}^2 \Big) \|e_h\|_{\mathbb{L}^2}^2 \\
&\quad + c h^{2\varrho} \|u\|_{\mathbb{L}^2}^2 \|u\|_{X_{1}}^2 \|u\|_{X_{\varrho}}^2 + c h^{2\varrho} (\log(1+1/h))^2 \|u\|_{X_{\varrho+1}}^2.
\end{aligned}
\]
Combining the bounds for $ I_3^{(1)} $ and $ I_3^{(2)} $ yields the desired estimate \cref{eq:I3-bound}.

\medskip\noindent
\textbf{Step 4. (Estimate of $I_4$).} We establish the bound
\begin{equation}
  \label{eq:I4-bound}
  I_4 \leqslant
  \frac{1}{8} \|e_h\|_{\mathbb{H}_{D}^1}^2
  + \frac{1}{2} \sum_{n=1}^\infty \|L_{\zeta_n,h}e_h\|_{\mathbb{L}^2}^2
  + c h^{2\varrho} (\log(1+1/h))^2 \|u\|_{X_{\varrho+1}}^2.
\end{equation}
We begin with the decomposition, valid for any $ n \geqslant 1 $:
\[
L_{\zeta_n,h}u_h - \mathcal{J}_h L_{\zeta_n}u
= L_{\zeta_n,h}e_h
+ L_{\zeta_n,h}(\mathcal{J}_h - I)u
+ (L_{\zeta_n,h} - \mathcal{J}_h L_{\zeta_n})u.
\]
Applying the Cauchy--Schwarz and Young's inequalities gives, for any $ \epsilon \in (0,1) $,
\begin{equation}
  \label{eq:I4-bound-1}
  \begin{aligned}
    I_4 &\leqslant
    \left( \frac{1}{2} + \epsilon \right) \sum_{n=1}^\infty \|L_{\zeta_n,h}e_h\|_{\mathbb{L}^2}^2 \\
    &\quad {} + \frac{c}{\epsilon}\sum_{n=1}^\infty \|L_{\zeta_n,h}(\mathcal{J}_h -I)u\|_{\mathbb{L}^2}^2
    + \frac{c}{\epsilon}\sum_{n=1}^\infty \|(L_{\zeta_n,h} - \mathcal{J}_h L_{\zeta_n})u\|_{\mathbb{L}^2}^2.
  \end{aligned}
\end{equation}
We now bound the second and third terms.
Applying \cref{eq:L-zeta-n-h-bound}, \cref{eq:I-Jh-2} with $ \alpha = (1-\varrho)/2 $,
and the condition \cref{eq:zeta-sum-cond}, we obtain 
\begin{equation}
  \label{eq:zq-1}
  \sum_{n=1}^\infty \|L_{\zeta_n,h}(\mathcal{J}_h-I)u\|_{\mathbb{L}^2}^2
  \leqslant c h^{2\varrho} (\log(1+1/h))^2 \|u\|_{X_{\varrho+1}}^2.
\end{equation}
For the third term we have
\begin{align*}
  \sum_{n=1}^\infty \|(L_{\zeta_n,h} - \mathcal{J}_hL_{\zeta_n})u\|_{\mathbb{L}^2}^2
  & =
  \sum_{n=1}^\infty \|(\mathcal{P}_h - \mathcal{J}_h\mathcal{P})((\zeta_n\cdot\nabla)u)\|_{\mathbb{L}^2}^2 \\
  & \stackrel{\text{(i)}}{\leqslant} ch^{-2}\sum_{n=1}^\infty \|(\mathcal{P}_h - \mathcal{J}_h\mathcal{P})((\zeta_n\cdot\nabla)u)\|_{X_{-1,h}}^2 \\
  & \stackrel{\text{(ii)}}{\leqslant} c h^{2\varrho} \sum_{n=1}^\infty \|(\zeta_n\cdot\nabla)u\|_{\mathbb{H}_{D}^{\varrho}}^2 \\
  & \stackrel{\text{(iii)}}{\leqslant} c h^{2\varrho} \sum_{n=1}^\infty \|\zeta_n\|_{\mathbb{W}^{1,\infty}}^2 \|u\|_{X_{\varrho+1}}^2,
\end{align*}
where (i) used the inverse estimate \cref{eq:inverse} with $(\theta_1,\theta_2)=(-1,0)$,
(ii) used \cref{eq:JhP-Ph-1} with $\alpha = (1-\varrho)/2$,
and (iii) used \cref{eq:zetan-nabla-v-H} with $\theta=\varrho$
and the continuous embedding $X_{\varrho+1} \hookrightarrow \mathbb{H}^{\varrho+1}$
(see \cref{eq:HD-to-H,eq:X-to-HD} with $\theta=\varrho+1$).
By the condition \cref{eq:zeta-sum-cond}, it follows that
\begin{equation}
  \label{eq:zq-2}
  \sum_{n=1}^\infty \|(L_{\zeta_n,h} - \mathcal{J}_h L_{\zeta_n})u\|_{\mathbb{L}^2}^2
  \leqslant c h^{2\varrho} \|u\|_{X_{\varrho+1}}^2.
\end{equation}
Substituting \eqref{eq:zq-1} and \eqref{eq:zq-2} into \cref{eq:I4-bound-1} yields,
for any $ \epsilon \in (0,1) $,
\[
I_4 \leqslant
\left( \frac{1}{2} + \epsilon \right) \sum_{n=1}^\infty \|L_{\zeta_n,h}e_h\|_{\mathbb{L}^2}^2
+ \frac{c h^{2\varrho}}{\epsilon} (\log(1+1/h))^2 \|u\|_{X_{\varrho+1}}^2.
\]
By \cref{eq:zeta-sum-cond} and \cref{eq:L-zeta-n-h-bound},
we also have
\[
  \sum_{n=1}^\infty \|L_{\zeta_n,h}e_h\|_{\mathbb{L}^2}^2 \leqslant c\|e_h\|_{\mathbb{H}_{D}^1}^2.
\]
Inserting this into the previous inequality yields
\[
I_4 \leqslant c\epsilon \|e_h\|_{\mathbb{H}_{D}^1}^2
+ \frac{1}{2} \sum_{n=1}^\infty \|L_{\zeta_n,h}e_h\|_{\mathbb{L}^2}^2
+ \frac{c h^{2\varrho}}{\epsilon} (\log(1+1/h))^2 \|u\|_{X_{\varrho+1}}^2.
\]
The desired bound \cref{eq:I4-bound} now follows by choosing $\epsilon > 0$ sufficiently small.

\medskip\noindent
\textbf{Step 5. (Estimate of the Martingale $M$).} We establish the following bound for the martingale term $M$:
\begin{equation}
    \label{eq:M-bound}
\mathbb{E} \big[ \|M\|_{ C([0,T];\mathbb{R})} \big] \leqslant \epsilon \|e_h\|_{L^2(\Omega; C([0,T]; \mathbb{L}^2))}^2 + \frac{c h^{2\varrho}}{\epsilon} (\log(1+1/h))^2, \quad \forall \epsilon > 0.
\end{equation}
By the regularity property \cref{eq:u-basic} and the stability estimate
\cref{eq:Jh-stab} with $\theta = 0$, we have
$\mathcal{J}_h u \in L^\infty\bigl(\Omega; C([0,T]; \mathbb{L}^2)\bigr)$.
Since $e_h = u_h - \mathcal{J}_h u$,
the triangle inequality in conjunction with the stability estimate
\cref{eq:uh-stab} yields
\[
  e_h \in L^\infty\bigl(\Omega; C([0,T]; \mathbb{L}^2)\bigr).
\]
Since \cref{eq:L-zeta-n-h-adjoint} implies 
$ \langle L_{\zeta_n,h} e_h, \, e_h \rangle = 0$ for all $n \geqslant 1$,
the identity $e_h = u_h - \mathcal{J}_hu$ gives
\[
\langle L_{\zeta_n,h} u_h, \, e_h \rangle
= \langle L_{\zeta_n,h} \mathcal{J}_hu, \, e_h \rangle,
\quad \forall n \geqslant 1.
\]
Consequently, the integrand in the definition of $ M $ can be estimated as
\[
\begin{aligned}
 | \langle L_{\zeta_n,h} u_h - \mathcal{J}_h L_{\zeta_n} u, \, e_h \rangle |
&= | \langle (L_{\zeta_n,h} \mathcal{J}_h - \mathcal{J}_h L_{\zeta_n}) u, \, e_h \rangle | \\
&\leqslant \| (L_{\zeta_n,h} \mathcal{J}_h - \mathcal{J}_h L_{\zeta_n}) u \|_{\mathbb{L}^2} \|e_h\|_{\mathbb{L}^2}.
\end{aligned}
\]
Squaring both sides and summing over $ n $, we find  
\begin{align*}
 \sum_{n=1}^\infty \bigl|\big\langle L_{\zeta_n,h} u_h - \mathcal{J}_h L_{\zeta_n} u, \,  e_h \big\rangle\bigr|^2 
&\leqslant \sum_{n=1}^\infty \| (L_{\zeta_n,h} \mathcal{J}_h  - \mathcal{J}_h L_{\zeta_n}) u \|_{\mathbb{L}^2}^2 \|e_h\|_{\mathbb{L}^2}^2 \\
&\leqslant c h^{2\varrho} (\log(1+1/h))^2 \|u\|_{X_{\varrho+1}}^2 \|e_h\|_{\mathbb{L}^2}^2,
\end{align*}
where the last inequality follows from the decomposition 
\[ 
L_{\zeta_n,h}\mathcal{J}_h - \mathcal{J}_h L_{\zeta_n} = L_{\zeta_n,h}(\mathcal{J}_h - I) 
+ (L_{\zeta_n,h} - \mathcal{J}_hL_{\zeta_n}), \quad n \geqslant 1, 
\]
and the estimates from \cref{eq:zq-1} and \cref{eq:zq-2}.
Applying the Burkholder--Davis--Gundy inequality now yields  
\begin{align*}
\mathbb{E} \big[ \|M\|_{C([0,T];\mathbb{R})} \big]
& \leqslant c \mathbb{E} \left[ \bigg( \int_0^T \sum_{n=1}^\infty \big|\big\langle L_{\zeta_n,h} u_h(t) - \mathcal{J}_h L_{\zeta_n} u(t), \,  e_h(t) \big\rangle\big|^2 \mathrm{d}t \bigg)^{1/2} \right] \\
& \leqslant c h^{\varrho} \log(1+1/h) \mathbb{E} \left[ \|u\|_{L^2(0,T; X_{\varrho+1})} \|e_h\|_{ C([0,T]; \mathbb{L}^2)} \right].
\end{align*}
Finally, applying Young's inequality with an arbitrary $\epsilon >0 $ and using the regularity
property \eqref{eq:u-L2} for $u$
yields the desired bound \eqref{eq:M-bound}.

\medskip\noindent
\textbf{Step 6. (Conclusion).} By \cref{eq:Ph} with $(\theta_1,\theta_2) = (0,1/2)$,
\cref{eq:I-Jh-1} with
$\alpha = (1-\varrho)/2$ and $\beta = 1/2$,
and the assumption $u_0 \in L_{\mathcal{F}_0}^2(\Omega; X_{1/2})$,
we have the initial error bound: 
\[
\|(\mathcal{P}_h - \mathcal{J}_h)u_0\|_{\mathbb{L}^2} \leqslant ch^{1/2} \|u_0\|_{X_{1/2}},
\quad \mathbb{P}\text{-a.s.}
\]  
Substituting this bound, together with the estimates for $ I_1, \dots, I_4 $ from \cref{eq:I1-bound,eq:I2-bound,eq:I3-bound,eq:I4-bound}, into the identity \cref{eq:eh-sde}, we deduce that,
$\mathbb{P}$-a.s., for all $ t \in [0,T] $,
\begin{align*}
  \|e_h(t)\|_{\mathbb{L}^2}^2 + \frac{1}{4} \int_0^t \|e_h(s)\|_{\mathbb{H}_{D}^1}^2 \,\mathrm{d}s
  &\leqslant ch \|u_0\|_{X_{1/2}}^2
  + c \int_0^t B_1(s) \|e_h(s)\|_{\mathbb{L}^2}^2 \,\mathrm{d}s \\
  & \quad + c h^{2\varrho} (\log(1+1/h))^2 \int_0^T B_2(s) \mathrm{d}s
  + \|M\|_{C([0,T];\mathbb{R})},
\end{align*}
where
\begin{align*}
  B_1(s) &:= \|u(s)\|_{\mathbb{L}^2}^2 \|u(s)\|_{X_{1}}^2 + \|u_h(s)\|_{\mathbb{L}^2}^2 \|u_h(s)\|_{\mathbb{H}_{D}^1}^2, \\
  B_2(s) &:= \|u(s)\|_{X_{\varrho+1}}^2 + \big(1 + \|u(s)\|_{\mathbb{L}^2}^2\big) \|u(s)\|_{X_{\varrho}}^2 \|u(s)\|_{X_{1}}^2.
\end{align*}
Define $ z(t) := \sup_{r \in [0,t]} \|e_h(r)\|_{\mathbb{L}^2}^2 $.
Then, $\mathbb{P}$-a.s.~for all $t \in [0,T]$, 
\begin{align*}
  z(t) + \frac{1}{4} \int_0^t \|e_h(s)\|_{\mathbb{H}_{D}^1}^2 \,\mathrm{d}s 
  &\leqslant ch\|u_0\|_{X_{1/2}}^2
  + c \int_0^t B_1(s) z(s) \,\mathrm{d}s \\
  & + c h^{2\varrho}(\log(1+1/h))^2 \int_0^T B_2(s) \mathrm{d}s 
   + 2\|M\|_{C([0,T];\mathbb{R})}.
\end{align*}
Applying Gronwall's inequality then yields, $\mathbb{P}$-a.s.,
\begin{align*}
  & z(T) + \frac{1}{4} \int_0^T \|e_h(s)\|_{\mathbb{H}_{D}^1}^2 \,\mathrm{d}s \\
  \leqslant{} & c \exp\Bigl( c \int_0^T B_1(s) \, \mathrm{d}s \Bigr) 
  \times \Bigl[
    h \|u_0\|_{X_{1/2}}^2 + \|M\|_{C([0,T];\mathbb{R})} 
    + h^{2\varrho}(\log(1+1/h))^2 \int_0^T B_2(s) \, \mathrm{d}s
  \Bigr] \\
  \leqslant{} & c \exp\left( c\|u\|_{C([0,T];\mathbb{L}^2)}^2 \|u\|_{L^2(0,T;X_{1})}^2 + c\|u_h\|_{C([0,T];\mathbb{L}^2)}^2 \|u_h\|_{L^2(0,T;\mathbb{H}_{D}^1)}^2 \right) \\
  & \, \times \bigg[ h\|u_0\|_{X_{1/2}}^2 + \|M\|_{C([0,T];\mathbb{R})}
    + h^{2\varrho}(\log(1+1/h))^2 \\
    & \qquad \times  \left( \|u\|_{L^2(0,T;X_{\varrho+1})}^2 + \big(1 + \|u\|_{C([0,T];\mathbb{L}^2)}^2\big) \|u\|_{C([0,T];X_{\varrho})}^2 \|u\|_{L^2(0,T;X_{1})}^2 \right)
  \bigg].
\end{align*}
By the regularity of $u$ in \cref{eq:u-basic} and the stability of $ u_h $ in \cref{thm:uh-stab},
it follows that, $\mathbb{P}$-a.s.,
\begin{align*}
   z(T) + \frac{1}{4}\int_0^T \|e_h(s)\|_{\mathbb{H}_{D}^1}^2 \,\mathrm{d}s 
   & \leqslant c \Bigl[ h\|u_0\|_{X_{1/2}}^2 
   + \|M\|_{C([0,T];\mathbb{R})}  \\
   & + h^{2\varrho}(\log(1+1/h))^2 \left( \|u\|_{L^2(0,T;X_{\varrho+1})}^2 + \|u\|_{C([0,T];X_{\varrho})}^2 \right) \Bigr].
\end{align*}
Taking expectations on both sides and invoking the initial condition
$u_0 \in L_{\mathcal{F}_0}^2(\Omega;X_{1/2})$ and the regularity estimates
\cref{eq:u-L2,eq:u-C}, we arrive at
\[
  \mathbb{E} \, z(T) + \|e_h\|_{L^2(\Omega \times [0,T]; \mathbb{H}_{D}^1)}^2
  \leqslant c \left( h^{2\varrho}(\log(1+1/h))^2 + \mathbb{E}\big[ \|M\|_{C([0,T];\mathbb{R})} \big] \right).
\]
Since $z(T) = \|e_h\|_{C([0,T];\mathbb{L}^2)}^2$, the preceding inequality reads
\[
  \|e_h\|_{L^2(\Omega; C([0,T];\mathbb{L}^2))}^2 + \|e_h\|_{L^2(\Omega \times [0,T];\mathbb{H}_{D}^1)}^2
  \leqslant c \left( h^{2\varrho}(\log(1+1/h))^2 + \mathbb{E}\big[\|M\|_{C([0,T];\mathbb{R})}\big] \right).
\]
Inserting the bound on $\mathbb{E}[\|M\|_{C([0,T];\mathbb{R})}]$ from
\cref{eq:M-bound} and absorbing the term
$\|e_h\|_{L^2(\Omega;C([0,T];\mathbb{L}^2))}^2$ into the left-hand side by
selecting $\epsilon > 0$ sufficiently small, and subsequently taking square
roots, we obtain the error estimate
\begin{equation}
  \label{eq:eh-bound}
  \|e_h\|_{L^2(\Omega; C([0,T];\mathbb{L}^2))} + \|e_h\|_{L^2(\Omega \times [0,T];\mathbb{H}_{D}^1)} \leqslant c h^{\varrho}\log(1+1/h).
\end{equation}
To bound the total error $u-u_h$, we decompose it as $u-u_h = (u-\mathcal{J}_h u) - e_h$.
By the triangle inequality together
with \eqref{eq:eh-bound}, we obtain
\begin{equation}
  \label{eq:u-uh-pre}
\begin{aligned}
  & \|u - u_h\|_{L^2(\Omega; C([0,T];\mathbb{L}^2))}
  + \|u - u_h\|_{L^2(\Omega \times [0,T];\mathbb{H}_{D}^1)} \\
  \leqslant{}& \|u - \mathcal{J}_h u\|_{L^2(\Omega; C([0,T];\mathbb{L}^2))}
  + \|u - \mathcal{J}_h u\|_{L^2(\Omega \times [0,T];\mathbb{H}_{D}^1)}
  + c\, h^{\varrho}\log(1+1/h).
\end{aligned}
\end{equation}
For the interpolation error $u - \mathcal{J}_h u$, invoking \cref{eq:I-Jh-1}
with $\alpha = (1-\varrho)/2$ and $\beta = \varrho$, in view of the regularity
property \cref{eq:u-C}, yields
\[
  \|u - \mathcal{J}_h u\|_{L^2(\Omega;C([0,T];\mathbb{L}^2))}
  \leqslant c\, h^{\varrho},
\]
while \cref{eq:I-Jh-2} with $\alpha = (1-\varrho)/2$, in view of the regularity
property \cref{eq:u-L2}, gives
\[
  \|u - \mathcal{J}_h u\|_{L^2(\Omega \times [0,T];\mathbb{H}_{D}^1)}
  \leqslant c\, h^{\varrho}\log(1+1/h).
\]
Inserting these two estimates into \cref{eq:u-uh-pre} yields the desired
error estimate \cref{eq:u-uh}. This completes the proof of
\cref{thm:u-uh}.

\hfill$\blacksquare$

\section{Conclusions} \label{sec:conclusions}
This paper has provided a rigorous numerical analysis of the two-dimensional incompressible stochastic Navier--Stokes
equations with transport noise under no-slip boundary conditions. We first investigated the regularity of the
solution and then established a strong convergence rate of $h^\varrho\log(1+1/h)$,
where $0 < \varrho < 1/2$, for a spatial semidiscretization based on the MINI finite element method.

To the best of our knowledge, this work constitutes the first rigorous convergence analysis for
finite element spatial semidiscretizations of the two-dimensional stochastic Navier--Stokes equations
with transport noise.
Our analysis relies on a smallness condition on the transport noise vector fields:
$\sum_{n=1}^\infty \|\zeta_n\|_{\mathbb{W}^{1,\infty}}^2 $ is sufficiently small.
Removing this restriction remains a fundamental open question in the numerical analysis of
stochastic fluid models. Key challenges for future work include extending these results to
fully discrete schemes and developing efficient algorithms to mitigate the substantial
computational cost of simulations.


\begin{thebibliography}{38}
\ifx \bisbn   \undefined \def \bisbn  #1{ISBN #1}\fi
\ifx \binits  \undefined \def \binits#1{#1}\fi
\ifx \bauthor  \undefined \def \bauthor#1{#1}\fi
\ifx \batitle  \undefined \def \batitle#1{#1}\fi
\ifx \bjtitle  \undefined \def \bjtitle#1{#1}\fi
\ifx \bvolume  \undefined \def \bvolume#1{\textbf{#1}}\fi
\ifx \byear  \undefined \def \byear#1{#1}\fi
\ifx \bissue  \undefined \def \bissue#1{#1}\fi
\ifx \bfpage  \undefined \def \bfpage#1{#1}\fi
\ifx \blpage  \undefined \def \blpage #1{#1}\fi
\ifx \burl  \undefined \def \burl#1{\textsf{#1}}\fi
\ifx \doiurl  \undefined \def \doiurl#1{\url{https://doi.org/#1}}\fi
\ifx \betal  \undefined \def \betal{\textit{et al.}}\fi
\ifx \binstitute  \undefined \def \binstitute#1{#1}\fi
\ifx \binstitutionaled  \undefined \def \binstitutionaled#1{#1}\fi
\ifx \bctitle  \undefined \def \bctitle#1{#1}\fi
\ifx \beditor  \undefined \def \beditor#1{#1}\fi
\ifx \bpublisher  \undefined \def \bpublisher#1{#1}\fi
\ifx \bbtitle  \undefined \def \bbtitle#1{#1}\fi
\ifx \bedition  \undefined \def \bedition#1{#1}\fi
\ifx \bseriesno  \undefined \def \bseriesno#1{#1}\fi
\ifx \blocation  \undefined \def \blocation#1{#1}\fi
\ifx \bsertitle  \undefined \def \bsertitle#1{#1}\fi
\ifx \bsnm \undefined \def \bsnm#1{#1}\fi
\ifx \bsuffix \undefined \def \bsuffix#1{#1}\fi
\ifx \bparticle \undefined \def \bparticle#1{#1}\fi
\ifx \barticle \undefined \def \barticle#1{#1}\fi
\bibcommenthead
\ifx \bconfdate \undefined \def \bconfdate #1{#1}\fi
\ifx \botherref \undefined \def \botherref #1{#1}\fi
\ifx \url \undefined \def \url#1{\textsf{#1}}\fi
\ifx \bchapter \undefined \def \bchapter#1{#1}\fi
\ifx \bbook \undefined \def \bbook#1{#1}\fi
\ifx \bcomment \undefined \def \bcomment#1{#1}\fi
\ifx \oauthor \undefined \def \oauthor#1{#1}\fi
\ifx \citeauthoryear \undefined \def \citeauthoryear#1{#1}\fi
\ifx \endbibitem  \undefined \def \endbibitem {}\fi
\ifx \bconflocation  \undefined \def \bconflocation#1{#1}\fi
\ifx \arxivurl  \undefined \def \arxivurl#1{\textsf{#1}}\fi
\csname PreBibitemsHook\endcsname

\bibitem[\protect\citeauthoryear{Agresti and Veraar}{2024}]{Agresti2024}
\begin{barticle}
\bauthor{\bsnm{Agresti}, \binits{A.}},
\bauthor{\bsnm{Veraar}, \binits{M.}}:
\batitle{Stochastic {Navier--Stokes} equations for turbulent flows in critical spaces}.
\bjtitle{Commun. Math. Phys.}
\bvolume{405}(\bissue{2}),
\bfpage{43}
(\byear{2024})
\end{barticle}
\endbibitem

\bibitem[\protect\citeauthoryear{Bessaih et~al.}{2014}]{Bessaih2014}
\begin{barticle}
\bauthor{\bsnm{Bessaih}, \binits{H.}},
\bauthor{\bsnm{Brze\'{z}niak}, \binits{Z.}},
\bauthor{\bsnm{Millet}, \binits{A.}}:
\batitle{{Splitting up method for the 2D stochastic Navier--Stokes equations}}.
\bjtitle{Stoch. Partial Differ. Equ. Anal. Comput.}
\bvolume{2}(\bissue{4}),
\bfpage{433}--\blpage{470}
(\byear{2014})
\end{barticle}
\endbibitem

\bibitem[\protect\citeauthoryear{Brze\'{z}niak et~al.}{1991}]{Brzezniak1991}
\begin{barticle}
\bauthor{\bsnm{Brze\'{z}niak}, \binits{Z.}},
\bauthor{\bsnm{Capi\'{n}ski}, \binits{M.}},
\bauthor{\bsnm{Flandoli}, \binits{F.}}:
\batitle{{Stochastic partial differential equations and turbulence}}.
\bjtitle{Math. Models Methods Appl. Sci.}
\bvolume{1}(\bissue{1}),
\bfpage{41}--\blpage{59}
(\byear{1991})
\end{barticle}
\endbibitem

\bibitem[\protect\citeauthoryear{Brze\'{z}niak et~al.}{2013}]{Prohl2013}
\begin{barticle}
\bauthor{\bsnm{Brze\'{z}niak}, \binits{Z.}},
\bauthor{\bsnm{Carelli}, \binits{E.}},
\bauthor{\bsnm{Prohl}, \binits{A.}}:
\batitle{{Finite-element-based discretizations of the incompressible {N}avier–{S}tokes equations with multiplicative random forcing}}.
\bjtitle{IMA J. Numer. Anal.}
\bvolume{33}(\bissue{3}),
\bfpage{771}--\blpage{824}
(\byear{2013})
\end{barticle}
\endbibitem

\bibitem[\protect\citeauthoryear{Breit and Dodgson}{2021}]{Breit2021}
\begin{barticle}
\bauthor{\bsnm{Breit}, \binits{D.}},
\bauthor{\bsnm{Dodgson}, \binits{A.}}:
\batitle{{Convergence rates for the numerical approximation of the {2D} stochastic {N}avier-{S}tokes equations}}.
\bjtitle{Numer. Math.}
\bvolume{147}(\bissue{3}),
\bfpage{553}--\blpage{578}
(\byear{2021})
\end{barticle}
\endbibitem

\bibitem[\protect\citeauthoryear{Boyer and Fabrie}{2013}]{Boyer2013}
\begin{bbook}
\bauthor{\bsnm{Boyer}, \binits{F.}},
\bauthor{\bsnm{Fabrie}, \binits{P.}}:
\bbtitle{Mathematical Tools for the Study of the Incompressible Navier--Stokes Equations and Related Models}.
\bpublisher{Springer},
\blocation{New York}
(\byear{2013})
\end{bbook}
\endbibitem

\bibitem[\protect\citeauthoryear{Bessaih and Millet}{2019}]{Bessaih2019}
\begin{barticle}
\bauthor{\bsnm{Bessaih}, \binits{H.}},
\bauthor{\bsnm{Millet}, \binits{A.}}:
\batitle{{Strong $L^2$ convergence of time numerical schemes for the stochastic two-dimensional Navier--Stokes equations}}.
\bjtitle{IMA J. Numer. Anal.}
\bvolume{39}(\bissue{4}),
\bfpage{2135}--\blpage{2167}
(\byear{2019})
\end{barticle}
\endbibitem

\bibitem[\protect\citeauthoryear{Bessaih and Millet}{2022}]{Bessaih2022}
\begin{barticle}
\bauthor{\bsnm{Bessaih}, \binits{H.}},
\bauthor{\bsnm{Millet}, \binits{A.}}:
\batitle{Strong rates of convergence of space-time discretization schemes for the 2{D} {Navier--Stokes} equations with additive noise}.
\bjtitle{Stoch. Dyn.}
\bvolume{22}(\bissue{2}),
\bfpage{2240005}
(\byear{2022})
\end{barticle}
\endbibitem

\bibitem[\protect\citeauthoryear{Breit et~al.}{2026}]{Breit2026}
\begin{barticle}
\bauthor{\bsnm{Breit}, \binits{D.}},
\bauthor{\bsnm{Moyo}, \binits{T.C.}},
\bauthor{\bsnm{Prohl}, \binits{A.}},
\bauthor{\bsnm{Wichmann}, \binits{J.}}:
\batitle{{Mean square temporal error estimates for the two-dimensional stochastic Navier--Stokes equations with transport noise}}.
\bjtitle{IMA J. Numer. Anal.}
\bvolume{46}(\bissue{3}),
\bfpage{1898}--\blpage{1938}
(\byear{2026})
\end{barticle}
\endbibitem

\bibitem[\protect\citeauthoryear{Breit and Prohl}{2023}]{Breit2023IMA}
\begin{barticle}
\bauthor{\bsnm{Breit}, \binits{D.}},
\bauthor{\bsnm{Prohl}, \binits{A.}}:
\batitle{{Numerical analysis of two-dimensional Navier–Stokes equations with additive stochastic forcing}}.
\bjtitle{IMA J. Numer. Anal.}
\bvolume{43}(\bissue{3}),
\bfpage{1391}--\blpage{1421}
(\byear{2023})
\end{barticle}
\endbibitem

\bibitem[\protect\citeauthoryear{Breit and Prohl}{2024}]{Breit2024}
\begin{barticle}
\bauthor{\bsnm{Breit}, \binits{D.}},
\bauthor{\bsnm{Prohl}, \binits{A.}}:
\batitle{{Error analysis for 2D stochastic Navier--Stokes equations in bounded domains with Dirichlet data}}.
\bjtitle{Found. Comput. Math.}
\bvolume{24}(\bissue{5}),
\bfpage{1643}--\blpage{1672}
(\byear{2024})
\end{barticle}
\endbibitem

\bibitem[\protect\citeauthoryear{Brenner and Scott}{2008}]{Brenner2008}
\begin{bbook}
\bauthor{\bsnm{Brenner}, \binits{S.C.}},
\bauthor{\bsnm{Scott}, \binits{L.R.}}:
\bbtitle{{The Mathematical Theory of Finite Element Methods}},
\bedition{3}rd edn.
\bpublisher{Springer},
\blocation{New York}
(\byear{2008})
\end{bbook}
\endbibitem

\bibitem[\protect\citeauthoryear{Bensoussan and Temam}{1973}]{Bensoussan1973}
\begin{barticle}
\bauthor{\bsnm{Bensoussan}, \binits{A.}},
\bauthor{\bsnm{Temam}, \binits{R.}}:
\batitle{{\'{E}}quations stochastiques du type navier--stokes}.
\bjtitle{J. Funct. Anal.}
\bvolume{13},
\bfpage{195}--\blpage{222}
(\byear{1973})
\end{barticle}
\endbibitem

\bibitem[\protect\citeauthoryear{Capi\'{n}ski and Cutland}{1991}]{Capinski1991}
\begin{barticle}
\bauthor{\bsnm{Capi\'{n}ski}, \binits{M.}},
\bauthor{\bsnm{Cutland}, \binits{N.J.}}:
\batitle{{Stochastic Navier--Stokes equations}}.
\bjtitle{Acta Appl. Math.}
\bvolume{25}(\bissue{1}),
\bfpage{59}--\blpage{85}
(\byear{1991})
\end{barticle}
\endbibitem

\bibitem[\protect\citeauthoryear{Carelli and Prohl}{2012}]{Prohl2012}
\begin{barticle}
\bauthor{\bsnm{Carelli}, \binits{E.}},
\bauthor{\bsnm{Prohl}, \binits{A.}}:
\batitle{{Rates of convergence for discretizations of the stochastic incompressible Navier--Stokes equations}}.
\bjtitle{SIAM J. Numer. Anal.}
\bvolume{50}(\bissue{5}),
\bfpage{2467}--\blpage{2496}
(\byear{2012})
\end{barticle}
\endbibitem

\bibitem[\protect\citeauthoryear{Dauge}{1989}]{Dauge1989}
\begin{barticle}
\bauthor{\bsnm{Dauge}, \binits{M.}}:
\batitle{{Stationary Stokes and Navier--Stokes systems on two- or three-dimensional domains with corners. Part I: linearized equations}}.
\bjtitle{SIAM J. Math. Anal.}
\bvolume{20}(\bissue{1}),
\bfpage{74}--\blpage{97}
(\byear{1989})
\end{barticle}
\endbibitem

\bibitem[\protect\citeauthoryear{Doghman}{2024}]{Doghman2024}
\begin{barticle}
\bauthor{\bsnm{Doghman}, \binits{J.}}:
\batitle{Numerical approximation of the stochastic {Navier--Stokes} equations through artificial compressibility}.
\bjtitle{Calcolo}
\bvolume{61}(\bissue{2}),
\bfpage{23}
(\byear{2024})
\end{barticle}
\endbibitem

\bibitem[\protect\citeauthoryear{Flandoli and Gatarek}{1995}]{Flandoli1995martingale}
\begin{barticle}
\bauthor{\bsnm{Flandoli}, \binits{F.}},
\bauthor{\bsnm{Gatarek}, \binits{D.}}:
\batitle{{Martingale and stationary solutions for stochastic Navier--Stokes equations}}.
\bjtitle{Probab. Theory Related Fields}
\bvolume{102}(\bissue{3}),
\bfpage{367}--\blpage{391}
(\byear{1995})
\end{barticle}
\endbibitem

\bibitem[\protect\citeauthoryear{Flandoli and Luo}{2021}]{Flandoli2021high}
\begin{barticle}
\bauthor{\bsnm{Flandoli}, \binits{F.}},
\bauthor{\bsnm{Luo}, \binits{D.}}:
\batitle{High mode transport noise improves vorticity blow-up control in 3d {Navier--Stokes} equations}.
\bjtitle{Probab. Theory Related Fields}
\bvolume{180}(\bissue{1--2}),
\bfpage{309}--\blpage{363}
(\byear{2021})
\end{barticle}
\endbibitem

\bibitem[\protect\citeauthoryear{Flandoli and Luongo}{2023}]{Flandoli2023}
\begin{bbook}
\bauthor{\bsnm{Flandoli}, \binits{F.}},
\bauthor{\bsnm{Luongo}, \binits{E.}}:
\bbtitle{Stochastic Partial Differential Equations in Fluid Mechanics}.
\bpublisher{Springer},
\blocation{Singapore}
(\byear{2023})
\end{bbook}
\endbibitem

\bibitem[\protect\citeauthoryear{Feng and Vo}{2025}]{feng2025full}
\begin{botherref}
\oauthor{\bsnm{Feng}, \binits{X.}},
\oauthor{\bsnm{Vo}, \binits{L.}}:
Full moment error estimates in strong norms for numerical approximations of stochastic navier--stokes equations with multiplicative noise. {Part~I}: Time discretization.
arXiv:2510.08291
(2025)
\end{botherref}
\endbibitem

\bibitem[\protect\citeauthoryear{Giga and Miyakawa}{1985}]{Giga1985B}
\begin{barticle}
\bauthor{\bsnm{Giga}, \binits{Y.}},
\bauthor{\bsnm{Miyakawa}, \binits{T.}}:
\batitle{{Solutions in $L_r$ of the Navier--Stokes initial value problem}}.
\bjtitle{Arch. Ration. Mech. Anal.}
\bvolume{89}(\bissue{3}),
\bfpage{267}--\blpage{281}
(\byear{1985})
\end{barticle}
\endbibitem

\bibitem[\protect\citeauthoryear{Girault and Raviart}{1986}]{Girault1986}
\begin{bbook}
\bauthor{\bsnm{Girault}, \binits{V.}},
\bauthor{\bsnm{Raviart}, \binits{P.-A.}}:
\bbtitle{{Finite Element Methods for Navier--Stokes Equations}}.
\bpublisher{Springer},
\blocation{Berlin}
(\byear{1986})
\end{bbook}
\endbibitem

\bibitem[\protect\citeauthoryear{Holm}{2015}]{Holm2015variational}
\begin{barticle}
\bauthor{\bsnm{Holm}, \binits{D.D.}}:
\batitle{{Variational principles for stochastic fluid dynamics}}.
\bjtitle{Proc. R. Soc. A}
\bvolume{471}(\bissue{2176}),
\bfpage{20140963}
(\byear{2015})
\end{barticle}
\endbibitem

\bibitem[\protect\citeauthoryear{Hausenblas and Randrianasolo}{2019}]{Hausenblas2019}
\begin{barticle}
\bauthor{\bsnm{Hausenblas}, \binits{E.}},
\bauthor{\bsnm{Randrianasolo}, \binits{T.}}:
\batitle{{Time-discretization of stochastic 2-D Navier--Stokes equations with a penalty-projection method}}.
\bjtitle{Numer. Math.}
\bvolume{143}(\bissue{2}),
\bfpage{339}--\blpage{378}
(\byear{2019})
\end{barticle}
\endbibitem

\bibitem[\protect\citeauthoryear{Hyt{\"o}nen et~al.}{2016}]{Hytonen2016}
\begin{bbook}
\bauthor{\bsnm{Hyt{\"o}nen}, \binits{T.}},
\bauthor{\bsnm{van Neerven}, \binits{J.}},
\bauthor{\bsnm{Veraar}, \binits{M.}},
\bauthor{\bsnm{Weis}, \binits{L.}}:
\bbtitle{{Analysis in Banach Spaces}: Volume {I}: Martingales and Littlewood--Paley Theory}.
\bpublisher{Springer},
\blocation{Cham}
(\byear{2016})
\end{bbook}
\endbibitem

\bibitem[\protect\citeauthoryear{Liu and R\"ockner}{2010}]{LiuWei2010}
\begin{barticle}
\bauthor{\bsnm{Liu}, \binits{W.}},
\bauthor{\bsnm{R\"ockner}, \binits{M.}}:
\batitle{{SPDE in Hilbert space with locally monotone coefficients}}.
\bjtitle{J. Funct. Anal.}
\bvolume{259}(\bissue{11}),
\bfpage{2902}--\blpage{2922}
(\byear{2010})
\end{barticle}
\endbibitem

\bibitem[\protect\citeauthoryear{Liu and R{\"o}ckner}{2015}]{Liu2015}
\begin{bbook}
\bauthor{\bsnm{Liu}, \binits{W.}},
\bauthor{\bsnm{R{\"o}ckner}, \binits{M.}}:
\bbtitle{Stochastic Partial Differential Equations: An Introduction}.
\bpublisher{Springer},
\blocation{Cham}
(\byear{2015})
\end{bbook}
\endbibitem

\bibitem[\protect\citeauthoryear{Mitrea and Monniaux}{2008}]{Mitrea2008}
\begin{barticle}
\bauthor{\bsnm{Mitrea}, \binits{M.}},
\bauthor{\bsnm{Monniaux}, \binits{S.}}:
\batitle{{The regularity of the Stokes operator and the Fujita--Kato approach to the Navier--Stokes initial value problem in Lipschitz domains}}.
\bjtitle{J. Funct. Anal.}
\bvolume{254}(\bissue{6}),
\bfpage{1522}--\blpage{1574}
(\byear{2008})
\end{barticle}
\endbibitem

\bibitem[\protect\citeauthoryear{Mikulevicius and Rozovskii}{2005}]{Mikulevicius2005}
\begin{barticle}
\bauthor{\bsnm{Mikulevicius}, \binits{R.}},
\bauthor{\bsnm{Rozovskii}, \binits{B.L.}}:
\batitle{{Global $L_2$-solutions of stochastic Navier--Stokes equations}}.
\bjtitle{Ann. Probab.}
\bvolume{33}(\bissue{1}),
\bfpage{137}--\blpage{176}
(\byear{2005})
\end{barticle}
\endbibitem

\bibitem[\protect\citeauthoryear{van Neerven et~al.}{2012}]{Neerven2012b}
\begin{barticle}
\bauthor{\bsnm{van Neerven}, \binits{J.}},
\bauthor{\bsnm{Veraar}, \binits{M.}},
\bauthor{\bsnm{Weis}, \binits{L.}}:
\batitle{{Maximal $L^p$-regularity for stochastic evolution equations}}.
\bjtitle{SIAM J. Math. Anal.}
\bvolume{44}(\bissue{3}),
\bfpage{1372}--\blpage{1414}
(\byear{2012})
\end{barticle}
\endbibitem



\bibitem[\protect\citeauthoryear{Ondrej\'at et~al.}{2023}]{Ondrejat2023}
\begin{barticle}
\bauthor{\bsnm{Ondrej\'at}, \binits{M.}},
\bauthor{\bsnm{Prohl}, \binits{A.}},
\bauthor{\bsnm{Walkington}, \binits{N.J.}}:
\batitle{{Numerical approximation of nonlinear SPDE's}}.
\bjtitle{Stoch. Partial Differ. Equ. Anal. Comput.}
\bvolume{11}(\bissue{4}),
\bfpage{1553}--\blpage{1634}
(\byear{2023})
\end{barticle}
\endbibitem

\bibitem[\protect\citeauthoryear{Pardoux and R{\u a}\c{s}canu}{2014}]{Pardoux2014}
\begin{bbook}
\bauthor{\bsnm{Pardoux}, \binits{E.}},
\bauthor{\bsnm{R{\u a}\c{s}canu}, \binits{A.}}:
\bbtitle{Stochastic Differential Equations, Backward SDEs, Partial Differential Equations}.
\bpublisher{Springer},
\blocation{Cham}
(\byear{2014})
\end{bbook}
\endbibitem

\bibitem[\protect\citeauthoryear{Runst and Sickel}{1996}]{Runst1996sobolev}
\begin{bbook}
\bauthor{\bsnm{Runst}, \binits{T.}},
\bauthor{\bsnm{Sickel}, \binits{W.}}:
\bbtitle{Sobolev Spaces of Fractional Order, Nemytskij Operators, and Nonlinear Partial Differential Equations}.
\bpublisher{Walter de Gruyter},
\blocation{Berlin}
(\byear{1996})
\end{bbook}
\endbibitem

\bibitem[\protect\citeauthoryear{R{\"o}ckner and Zhang}{2009}]{Rockner2009}
\begin{barticle}
\bauthor{\bsnm{R{\"o}ckner}, \binits{M.}},
\bauthor{\bsnm{Zhang}, \binits{X.}}:
\batitle{Stochastic tamed 3{D} {N}avier-{S}tokes equations: existence, uniqueness and ergodicity}.
\bjtitle{Probab. Theory Related Fields}
\bvolume{145}(\bissue{1--2}),
\bfpage{211}--\blpage{267}
(\byear{2009})
\end{barticle}
\endbibitem

\bibitem[\protect\citeauthoryear{Scheutzow}{2013}]{Scheutzow2013}
\begin{barticle}
\bauthor{\bsnm{Scheutzow}, \binits{M.}}:
\batitle{A stochastic {G}ronwall lemma}.
\bjtitle{Infin. Dimens. Anal. Quantum Probab. Related Top.}
\bvolume{16}(\bissue{2}),
\bfpage{1350019}
(\byear{2013})
\end{barticle}
\endbibitem

\bibitem[\protect\citeauthoryear{Sohr}{2012}]{Sohr2012navier}
\begin{bbook}
\bauthor{\bsnm{Sohr}, \binits{H.}}:
\bbtitle{The Navier-Stokes Equations: An Elementary Functional Analytic Approach}.
\bpublisher{Birkh{\"a}user},
\blocation{Basel}
(\byear{2012})
\end{bbook}
\endbibitem


\bibitem[\protect\citeauthoryear{Yagi}{2010}]{Yagi2010}
\begin{bbook}
\bauthor{\bsnm{Yagi}, \binits{A.}}:
\bbtitle{{Abstract Parabolic Evolution Equations and Their Applications}}.
\bpublisher{Springer},
\blocation{Berlin}
(\byear{2010})
\end{bbook}
\endbibitem

\end{thebibliography}


\end{document}